\documentclass[]{interact}

\usepackage{epstopdf}
\usepackage{booktabs}
\usepackage{caption}
\usepackage{subcaption}
\usepackage{enumitem}

\usepackage{url}

\usepackage[numbers,sort&compress]{natbib}
\bibpunct[, ]{[}{]}{,}{n}{,}{,}
\makeatletter
\def\NAT@def@citea{\def\@citea{\NAT@separator}}
\makeatother

\usepackage{xcolor}
\usepackage{amsmath,amsthm,amssymb,bm}
\usepackage{mathtools}
\usepackage{algorithm}
\usepackage{algpseudocode}

\usepackage[mathscr]{eucal}

\usepackage[official]{eurosym}

\usepackage{pgf,tikz}

\usetikzlibrary{calc}
\usetikzlibrary{shapes.geometric,arrows,shapes,chains,automata}

\def\N{{\mathbb{N}}}
\def\R{{\mathbb{R}}}

\def\d{{\mathcal{D}}}

\def\m{{\mathcal{M}}}
\def\M{{\mathscr{M}}}

\def\Ie{{\mathscr I_e}}

\def\tu{{\tilde{u}}}
\def\tx{{\tilde{x}}}
\def\ty{{\tilde{y}}}
\def\tJ{{\tilde{J}}}

\def\u{{\mathscr U}}

\def\uad{{\mathscr U_{\mathsf{ad}}}}
\def\Uad{{U_{\mathsf{ad}}}}
\def\x{{\mathscr X}}
\def\y{{\mathscr Y}}
\def\xe{{\mathscr X_e}}
\def\ye{{\mathscr Y_e}}
\def\xad{{\mathscr X_0}}
\def\nz{{n_{\mathsf z}}}
\def\nx{{n_{\mathsf x}}}
\def\ny{{n_{\mathsf y}}}
\def\nu{{n_{\mathsf u}}}
\def\ua{{u_{\mathsf a}}}
\def\ub{{u_{\mathsf b}}}
\def\dx{{\mathscr{D}_{\mathsf x}}}
\def\dy{{\mathscr{D}_{\mathsf y}}}
\def\dz{{\mathscr{D}_\mathsf z}}

\def\d{{\mathscr{D}}}

\def\tilx0{{\tilde{x}_0}}

\def\eps{{\varepsilon}}

\numberwithin{equation}{section}

\theoremstyle{plain}
\newtheorem{theorem}{Theorem}[section] 

\theoremstyle{definition}
\newtheorem{definition}[theorem]{Definition} 
\newtheorem{lemma}[theorem]{Lemma}

\newtheorem*{theorem*}{Theorem}

\newtheorem{examp}[theorem]{Example} 

\newtheorem*{examp*}{Example} 
\newenvironment{example*}
  {\renewcommand{\qedsymbol}{$\Diamond$}%
   \pushQED{\qed}\begin{examp*}}
  {\popQED\end{examp*}}
  
\newtheorem{remark}[theorem]{Remark}  
  
\newtheorem{assumption}[theorem]{Assumption}  

\newtheorem{con}[theorem]{Conjecture}  

\newtheorem{pro}[theorem]{Problem}  

\begin{document}


\title{Optimal Control of Dynamic District Heating Networks}

\author{
\name{Christian~J\"akle\textsuperscript{a}, Lena~Reichle\textsuperscript{a}, and Stefan~Volkwein\textsuperscript{a}\thanks{S.~Volkwein. Email: Stefan.Volkwein@uni-konstanz.de}}
\affil{\textsuperscript{a}Department of Mathematics and Statistics, University of Konstanz, 78457 Konstanz, Germany}
}

\maketitle

\begin{abstract}
	In the present paper an optimal control problem for a system of differential-algebraic equations (DAEs) is considered. This problem arises in the dynamic optimization of unsteady district heating networks. Based on the Carath\'eodory theory existence of a unique solution to the DAE system is proved using specific properties of the district heating network model. Moreover, it is shown that the optimal control problem possesses optimal solutions. For the numerical experiments different networks are considered including also data from a real district heating network.
\end{abstract}


\section{Introduction}
\label{sec:1}

The efficient use of energy, particularly renewable energy sources, plays an important role in today’s discussions. Therefore district heating and particular its optimization become more important in energy use, as it is flexible in the supply of different forms of energy. District heating is a system that transports heat energy via a network of pipelines from a central power plant to different consumers with different requests. For a long time, district heating has been considered as a static problem, due to the volatility and diversity in energies and the different requests of the consumers, this assumption becomes obsolete. Therefore, it has become more important to simulate and optimize the dynamic processes based on changes in supply/demand of the consumers and energy productions.

To realize an efficient use of district heating networks, we require (i) a suitable mathematical model of the network as well as (ii) fast and stable simulation and (iii) optimization techniques; cf. \cite{Krug19}. 
The mathematical model consists in a system of one-dimensional (1d) coupled hyperbolic partial differential equations (PDEs) to model water and heat transport. Moreover, the inclusion of a system of algebraic equations, representing the depot and households, introduces an additional level of non-linearity to the network. By applying principles of mass conservation and pressure continuity at each node, the model extends from individual pipes to the entire network.
It is explained in \cite{JaekleReichle} that the obtained mathematical model can be described by so-called differential algebraic equations (DAEs); cf., e.g., \cite{gerdts2011optimal,KM06}. Here we rely on  \cite{Borsche19,Krug19}.

The literature about the optimization of district heating networks is very limited. In \cite{pirouti2013energy}, the authors present an applied case study focusing on a particular district heating network located in South Wales. The study explores the specific context and dynamics of this network. On the other hand, \cite{rezaie2012district} offers a broader and more general discussion concerning the technological aspects and the potential benefits associated with district heating networks in general.
In \cite{schweiger2017district}, the authors adopt a "first-discretize-then-optimize" approach when dealing with the PDE-constrained problem at hand. This means they first discretize the underlying PDEs and then proceed with the optimization process.
To explore the connection between district heating networks and energy storage aspects, please see references \cite{colella2012numerical, hart2011pyomo, verda2011primary}, as well as the cited sources within those references for further details.
Studies on the design and expansion of networks also take into account stationary models of hot water flow. These models are explored in works such as \cite{BGV16, dorfner2014large, RS20, DMRS23}. 
In \cite{Borsche19}, the authors investigate the numerical simulation of district heating networks by employing a local time stepping method. On the other hand, the authors in \cite{rein2020optimal} as well as in \cite{RMDK21, rein2018parametric} delve into model order reduction techniques applied to the hyperbolic equations in district heating networks. These studies focus on enhancing the efficiency and accuracy of simulations in this context.
Lastly, in reference \cite{hauschild2020port}, the authors introduce and explore a port-Hamiltonian modeling approach specifically tailored for district heating networks. This approach offers a unique perspective and insights into the dynamics and control of such networks.

It turns out that the optimal control of dynamic district heating networks leads to an optimal control problem governed by differential algebraic equations (DAEs):
\begin{align}
    \nonumber
    &\min J(x,y,u) \coloneqq \int_{t_0}^{t_f} f_0(t,x(t),y(t),u(t))\, \mathrm dt\\
    \tag{\textbf P}
    \label{eq:DAE_Orig}
    &\hspace{1.0mm}\text{s.t. }(x, y, u) \text{ satisfies }u\in\uad\text{ and}\\
    &\nonumber\hspace{1mm}\mathrm{(1.1)}\hspace{4mm}\left\{\hspace{4mm} 
    \begin{aligned}
        \dot{x}(t)&=f(t,x(t),y(t),u(t)) && \text{f.a.a. } t\in(t_0, t_f],\quad x(t_0)=x_0,\\
        0 &= g(t, x(t), y(t), u(t)) && \text{f.a.a. }t\in(t_0, t_f],
    \end{aligned}
    \right.
\end{align}
where `f.a.a.' stands for `for almost all', $f_0 : [t_0,t_f]\times \R^\nx \times \R^\ny 
\times \R^\nu \to \R$, $f : [t_0,t_f] \times \R^\nx \times \R^\ny \times \R^\nu \to \R^\nx$, $g : [t_0,t_f] \times \R^\nx \times \R^\ny \times \R^\nu \to \R^\ny$, $x_0\in\xad \subset\R^\nx$ hold, with $\xad$ the set of admissible initial values and $\uad$ is a closed and convex set. For more details we refer the reader to Sections~\ref{sec:2} 
and \ref{sec:3}.

The new aspects of the present paper are as follows: (i) We prove existence of a unique solution pair $(x,y)$ to (1.1) in the extended sense applying the Carath\'eodory theory. Here we essentially utilize the specific structure of the DAE in the district heating application. (ii) Without assuming any growth rates we ensure the existence of optimal controls in $L^2(t_0,t_f;\R^\nu)$ for \eqref{eq:DAE_Orig}. (iii) In our numerical experiments we solve \eqref{eq:DAE_Orig} for a part of a real dynamic district heating network, where we get the data from the company Rechenzentrum f\"ur Versorgungsnetze Wehr GmbH\footnote{See \url{https://www.rzvn.de}}. Here we make use of an instantaneous optimal control approach to get consistent starting values for the optimization method.

The paper is organized in the following manner: In Section~\ref{sec:2} we prove existence of a unique solution to a semi-explicit DAE with differentiation index $1$ locally in time. The optimal control problem 
is studied in Section~\ref{sec:3}. The district heating network is formulated in the form 
\eqref{eq:DAE_Orig} in Section~\ref{sec:4}, whereas different numerical experiments 
are discussed in Section~\ref{sec:5}.

\section{Existence and uniqueness of a solution to the state equations}
\label{sec:2}

Suppose that $0\le t_0<t_f$ are fixed. The convex, closed and bounded set of admissible controls is given as
\begin{align*}
    \uad=\big\{u\in\u\,\big|\,u(t)\in\Uad\big\},
\end{align*}
where $\u$ stands for the Hilbert space $L^2(t_0,t_f;\R^\nu)$ endowed with the usual topology, $\Uad=\{\mathsf u\in\R^\nu\,|\,\ua\le\mathsf u\le\ub\text{ in }\R^\nu\}$ and $\ua\leq \ub\in\R^\nu$ are constants. Here `$\le$' is understood component-wise.

In a first step we study the state equations of \eqref{eq:DAE_Orig} -- i.e., the DAE -- and prove existence (local in time) of a unique solution. Throughout this section the control $u\in\uad$ is considered to be fixed. For fixed $u\in\uad$ let us introduce the two functions
\begin{align}\label{eq:modified_g}
    \left.
    \begin{aligned}
        \mathsf f(t,\mathsf x,\mathsf y)&=f(t,\mathsf x,\mathsf y,u(t))\hspace{2mm}\\
        \mathsf g(t,\mathsf x,\mathsf y)&=g(t,\mathsf x,\mathsf y,u(t))
    \end{aligned}
    \right\}\quad\text{f.a.a. }t\in[t_0,t_f]\text{ and for }(\mathsf x,\mathsf y)\in\R^\nx\times\R^\ny.
\end{align}
Utilizing 
Carath\'eodory's theory (cf. \cite[Appendix C, Theorem 54]{sontag1998}) we will show the well-posedness of the DAE
\begin{align}\label{eq:semiDAE}
	\begin{aligned}
	    \dot x(t)&=\mathsf f(t,x(t),y(t))&&\text{f.a.a. }t\in(t_0,t_f],\quad x(t_0)=x_0,\\
        0 &=\mathsf g(t,x(t),y(t))&&\text{f.a.a. } t \in [t_0,t_f].
	\end{aligned}
\end{align}
In addition, we also consider the ordinary differential equation (ODE) 
\begin{align}
    \label{eq:ivp_def}
    \dot z(t)&=\mathsf F(t,z(t))\text{ f.a.a. } t \in (t_0,t_f],\quad z(t_0)= z_0
\end{align}
with $z_0 \in \R^\nz$, $z(t) \in \R^\nz$ and $\mathsf F:[t_0,t_f]\times\R^\nz\to\R^\nz$.

\begin{definition}
    For $t_e\in(t_0,t_f]$ and $\Ie=[t_0,t_e]$ a function $\varphi:\Ie\to \R^n$ is \textit{absolutely continuous} 
    on $\Ie$ if for every $\eps > 0$ 
    there exists a $\delta > 0$ such that for every finite sequence 
    of pairwise disjoint subintervals $\{(a_k,b_k)\}_{k=1}^M$ with $a_k < b_k$, $(a_k,b_k)\subset\Ie$ and 
    $\sum_{k=1}^M |b_k-a_k| < \delta$ 
    it follows 
    \begin{align*}
        \sum_{k=1}^M {\|\varphi(b_k)-\varphi(a_k)\|}_2 < \eps,
    \end{align*}
    where $\|\cdot\|_2$ stands for the Euclidean norm. For $n\in\mathbb N$ we define the function space
    \begin{align*}
         AC(\Ie;\R^n):=\big\{\varphi:\Ie\to \R^n \,\big|\,\varphi \text{ is absolutely continuous on }\Ie\big\}.
    \end{align*}
\end{definition}

Absolutely continuous functions have the following property \cite[Appendix C.1, p.~471]{sontag1998}.

\begin{lemma}
    \label{lem:ac_derivative}
    Let $\varphi\in AC(\Ie;\R^n)$. Then $\varphi$ has a derivative $\varphi'$ almost everywhere. The derivative $\varphi'$ is Lebesgue integrable, and we have
    \begin{align*}
		\varphi(t)=\varphi(t_0) +\int_{t_0}^t \varphi'(s)\,\mathrm ds \quad\text{for all } t \in\Ie.
	\end{align*}
\end{lemma}

Due to Lemma~\ref{lem:ac_derivative} the norm
\begin{align*}
    {\|\varphi\|}_{AC(\Ie;\R^n)}=\max_{t\in\Ie}{\|\varphi(t)\|}_2+\int_{t_0}^{t_\mathsf e}{\|\varphi'(t)\|}_2\,\mathrm dt\quad\text{for }\varphi\in AC(\Ie;\R^n)
\end{align*}
is well-defined. Further, it is known that $AC([t_0,t_\mathsf e];\R^n)$ is a Banach space supplied with the norm $\|\cdot\|_{AC(\Ie;\R^n)}$; cf. \cite[Theorem~7.16]{Leo17}. Every $\varphi\in W^{1,1}(\Ie;\R^n)$ is equal almost everywhere to a function in $AC(\Ie;\R^n)$. We set
\begin{align*}
    \xe:=AC(\Ie;\R^\nx)\quad\text{and}\quad\ye:=L^2(\Ie;\R^\ny).
\end{align*}

\begin{definition}
    Let $t_e\in(t_0,t_f]$ and $\Ie=[t_0,t_e]$.
    \begin{enumerate}
        \item [1)] The function $z:\Ie\to \R^\nz$ is a \textit{solution in the extended sense} of the ODE \eqref{eq:ivp_def} if $z \in AC(\Ie;\R^\nz)$ holds and $z$ fulfills \eqref{eq:ivp_def} on $\Ie$.
        \item [2)] The function pair $(x,y):\Ie\to \R^\nx\times\R^\ny$ is a 
        \textit{solution in the extended sense} of the DAE \eqref{eq:semiDAE}, if we have $(x,y) \in \xe\times\ye$ and \eqref{eq:semiDAE} is satisfied on $\Ie$. 
  \end{enumerate}
\end{definition}

In the context of DAEs the notion of the differentiation index is very important. Therefore, we recall its definition next. Let us mention that in this work we suppose that \eqref{eq:semiDAE} has differentiation index $di=1$.

\begin{definition}
   Assume that $\dx\subset\R^\nx$ and $\dy\subset\R^\ny$ are given non-empty and open subsets and set $\d:=\dx\times\dy$. Let $t_e\in(t_0,t_f]$ be given and $\Ie=[t_0,t_e]$. We assume that for any $t\in\Ie$ the function $\mathsf g(t,\cdot\,,\cdot)$ is continuously differentiable in $\d$. Then, \eqref{eq:semiDAE} 
    has \textit{differentiation index} $di = 1$ on $\Ie\times \d$ if the Jacobian matrix
    \begin{align*}
        \partial_\mathsf y \mathsf g(t,\mathsf x,\mathsf y)\in\R^{\ny\times\ny}
    \end{align*}
    is regular for all $(\mathsf x,\mathsf y)\in\d$ and f.a.a. $t \in\Ie$. If \eqref{eq:semiDAE} has differentiation index $di = 1$, the initial value $x_t\in\dx$ is said to be \emph{consistent} at $t \in\Ie$ if there exists a $y_t\in\dy$ with $\mathsf g(t,x_t,y_t)=0$.
\end{definition}
We define the \emph{set of all consistent initial value pairs} at the starting time $t$ as
\begin{align*}
	\M_0(t):=\big\{\mathsf x\in\dx\,|\,\text{there exists a } \mathsf y \in \dy\text{ with }\mathsf g(t,\mathsf x,\mathsf y) = 0\}\quad\text{for }t \in\Ie.
\end{align*}
\begin{remark}
    \label{Rem-PL}
    \begin{enumerate}
        \item [1)] Based on Picard-Lindel\"of's theorem (cf., e.g., \cite[Chapter 9]{gerhardt06}) we already proved in \cite[Theorem 2]{JaekleReichle} that there exists a $t_e\in(t_0,t_f]$ such that \eqref{eq:semiDAE} has a unique, continuously differentiable solution pair $(x,y)$ on $\Ie=[t_0,t_e]$. Here, we have to assume that there are open subsets $\dx\subset\mathbb R^\nx$ and $\dy\subset\mathbb R^\ny$ so that $\mathsf f$ is continuous in $t$, uniformly Lipschitz-continuous in $(\mathsf x,\mathsf y)$ and that $\mathsf g$ is differentiable with a uniformly Lipschitz-continuous derivative on $\Ie\times \dx\times\dy$. However, in our optimal control problem these conditions are too strong. Due to \eqref{eq:modified_g} we cannot expect that both $\mathsf f$ and $\mathsf g$ are continuous with respect to $t$ for a given control function $u\in\uad$. Therefore, the assumptions we made in \cite{JaekleReichle} are usually no longer fulfilled.
        \item [2)] If the differentiation index $di$ is two, we need to handle $\partial_t \mathsf g(t,\mathsf x,\mathsf y)$. There are several methods to do that, one logical idea is to assume more regularity from the control, e.g., $u \in H^1(t_0,t_f;\R^\nu)$. Another idea is to require that the algebraic equation does neither depend on the control $u$ nor the algebraic variable $y$ and the control $u$ only occurs in the differential equation. Then we can write the algebraic equation in the form $\mathsf g(t,x(t)) = 0$; see \cite{gerdts2006local}, for instance. In this setting we can choose $u \in\uad$. After an index reduction, the results of our work can be applied. However, it must be taken into account that the solution of the reduced system still solves the original problem.\hfill$\Diamond$
    \end{enumerate}
\end{remark}

Next we start with the initial value problem \eqref{eq:ivp_def}. For given $\mathsf z_\circ\in\R^\nz$ and $\rho_\circ>0$ we denote by $B_{\rho_\circ}(\mathsf z_\circ)\subset\R^\nz$ the open ball centered at $\mathsf z_\circ$ with radius $\rho_\circ$. A proof of the next theorem can be found in \cite[Appendix C, Theorem 54]{sontag1998}, for instance.

\begin{theorem}[Carathéodory]
\label{thm:caratheo_ex}
	Let $\dz\subset\R^\nz$ be an open subset and the function $F:[t_0,t_f] \times \dz \to \R^n$ fulfill the following conditions: 
    \begin{enumerate}
		\item [\em 1)] $\mathsf F(\cdot\,,\mathsf z)$ is measurable on $[t_0,t_f]$ for all $\mathsf z \in \dz$;
		\item [\em 2)] $\mathsf F(t,\cdot)$ is continuous on $\dz$ f.a.a. $t \in [t_0,t_f]$;
		\item [\em 3)] for all $\mathsf z_\circ\in \dz$ there are $\rho_\circ>0$ and $\alpha\in L^1(t_0,t_f)$ such that $B_{\rho_\circ}(\mathsf z_\circ) \subset \dz$ and 
        \begin{align*}
            {\|\mathsf F(t,\mathsf z)-\mathsf F(t,\tilde{\mathsf z})\|}_2\le\alpha(t)\,{\|\mathsf z-\tilde{\mathsf z}\|}_2\quad\text{f.a.a. }t\in[t_0,\tilde t_f],\text{ for all }\mathsf z,\tilde{\mathsf z}\in B_{\rho_\circ}(\mathsf z_\circ);
        \end{align*}
        \item [\em 4)] for all $\mathsf z \in \dz$ there exists $\beta\in L^1(t_0,t_f)$ satisfying 
        \begin{align*}
            {\|\mathsf F(t,\mathsf z)\|}_2\leq \beta(t)\quad\text{f.a.a. }t\in[t_0,t_f].
        \end{align*}
	\end{enumerate}
	Then, for each initial condition $z_0 \in \dz$ there is a $t_e\in(t_0,t_f]$, so that \eqref{eq:ivp_def} possesses a unique solution $z$ in the 
	extended sense on $\Ie=[t_0,t_e]$ with $z(t)\in\dz$ f.a.a. $t\in\Ie$. 
\end{theorem}

\begin{remark}
    \label{rem:goal}
    \begin{enumerate}
        \item [1)] For $\rho_\circ>0$ let $\overline{B_{\rho_\circ}(z_0)}\subset\dz$ hold. We infer from Theorem~\ref{thm:caratheo_ex} that the solution $z$ is absolutely continuous in $\Ie$ and $z(t_0)\in B_{\rho_\circ}(z_0)$. Hence, we can assume (after possibly reducing $t_e>t_0$) that $z(t)\in B_{\rho_\circ}(z_0)$ holds for all $t\in\Ie$. This implies
        \begin{align*}
            {\|z\|}_{C(\Ie;\R^n)}&=\sup_{t\in\Ie}{\|z(t)\|}_2\le\sup_{t\in\Ie}\Big({\|z_0\|}_2+{\|z(t)-z_0\|}_2\Big)\le c_\mathsf{apr}
        \end{align*}
        with the a-priori bound $c_\mathsf{apr}=||z_0\|_2+\rho_0>0$.
        \item [2)] To prove existence results for the DAE based on Carathéodory's existence theorem one need to handle the algebraic equation 
        \begin{align}
            \label{eq:algebraic_equation}
		  \mathsf g(t,x(t),y(t)) = 0\quad\text{f.a.a. }t\in(t_0,t_f],
        \end{align}
        where $\mathsf g$ has been introduced in \eqref{eq:modified_g}. Here are different approaches possible:
        \begin{enumerate}
            \item [(a)] If the control function $u$ is sufficiently regular, the mapping $\mathsf g$ is continuously differentiable with respect to $(t,\mathsf x)$. In that case we can utilize the implicit function theorem to get a locally defined function $\mathsf y(t,\mathsf x)$ with
            \begin{align}
               \label{g_y}
                \mathsf g\big(t,x(t),\mathsf y(t,x(t))\big)=0\quad\text{f.a.a. }t\in\Ie.
            \end{align}
            As already mentioned in Remark~\ref{Rem-PL}-1), our control $u$ belongs to $\uad$ and consequently $\mathsf g$ is not continuously differentiable with respect to $t$.
            \item [(b)] For the case $u\in\uad$ we suppose that there exists a function $\mathsf y$ satisfying \eqref{g_y}. Thus, \eqref{eq:algebraic_equation} holds for $y(t)=\mathsf y(t,x(t))$. We prove that for our district heating application this assumption is fulfilled; cf. Lemma \ref{lem:representation}.
            \item [(c)] In \cite{roubivcek2002optimal} the authors request linear growth conditions for $\mathsf y=\mathsf y(t,\mathsf x)$. Based on them they prove the existence of solutions utilizing the existence theory of Filippov-Roxin \cite{filippov1962certain, roxin1962existence}. Since our applications do not satisfy linear growth conditions, we follow a different approach here.\hfill$\Diamond$
        \end{enumerate}
	\end{enumerate}
\end{remark}

\begin{assumption}\label{ass:representation_y}
	There exist an open set $\dx \subset \R^\nx$ and a function $\mathsf y\in L^2_{loc}((t_0,t_f)\times\dx;\R^\ny)$ satisfying
    \begin{enumerate}
		\item [1)] $\mathsf y(\cdot\,,\mathsf x) \in L^2(t_0,t_f;\R^\ny)$ for all $\mathsf x\in\dx$;
        \item [2)] $\mathsf g(t,\mathsf x,\mathsf y(t,\mathsf x)) =0$ f.a.a. $t\in[t_0,t_f]$ and for all $\mathsf x\in\dx$;
        \item [3)] $\mathsf y(t,\mathsf x)\in \dy$ f.a.a. $t\in[t_0,t_f]$ and for all $\mathsf x\in\dx$;
		\item [4)] $\mathsf f(\cdot\,,\mathsf x,\mathsf y(\cdot\,,\mathsf x))$ is measurable on $[t_0,t_f]$ for all $\mathsf x \in \dx$;
		\item [5)] $\mathsf y(t,\cdot)$ is continuous on $\dx$ f.a.a. $t\in[t_0,t_f]$;
		\item [6)] for all $\mathsf x_\circ\in \dx$ there are $\rho_\circ > 0$ and $\alpha_\mathsf y\in L^2(t_0,t_f)$ with 
        \begin{align*}
            {\|\mathsf y(t,\mathsf x)-\mathsf y(t,\tilde{\mathsf x})\|}_2 \leq \alpha_\mathsf y(t)\,{\|\mathsf x-\tilde{\mathsf x}\|}_2
        \end{align*}
        f.a.a. $t \in [t_0,t_f]$ and for all $\mathsf x,\tilde{\mathsf x} \in B_{\rho_\circ}(\mathsf x_\circ)$;
 		\item [7)] for all $\mathsf x \in \dx$ there exists $\beta\in L^1(t_0,t_f)$
        with 
        \begin{align*}
            {\|\mathsf f(t,\mathsf x,\mathsf y(t,\mathsf x))\|}_2 \leq \beta(t)
        \end{align*}
        f.a.a. $t \in [t_0,t_f]$.
	\end{enumerate}
\end{assumption}
\begin{remark}
    \begin{enumerate}
        \item [1)] In our work we assume the existence of the function $\mathsf y$ on $[t_0,t_f]\times \dx$. This can be verified for our specific application introduced in Section~\ref{sec:4}. It is also straightforward to only work locally. Then the choice of $\dx$ depends on the initial value $x_0\in\M_0(t_0)$. However, to simplify the presentation we suppose that we can work globally.
        \item [2)] Assumption~\ref{ass:representation_y}-3) can be ignored provided $\dy =\R^\ny$ holds. 
        \item [3)] If $\mathsf f(\cdot,\mathsf x,\mathsf y(\cdot\,,\mathsf x))$ is continuous on $[t_0,t_f]$ for all $\mathsf x \in \dx$, Assumption~\ref{ass:representation_y}-4) can be replaced by the condition that $\mathsf y(\cdot,\mathsf x)$ is measurable on $[t_0,t_f]$ for all $\mathsf x \in \dx$.\hfill$\Diamond$
    \end{enumerate}
\end{remark}
\begin{theorem}
    \label{thm:existence_index1}
	Let $\dx\in\mathbb R^\nx$ and $\dy\in\mathbb R^\ny$ be non-empty, open subsets and $\d:=\dx\times\dy$. In \eqref{eq:semiDAE} we consider $\mathsf f$ as a mapping from $[t_0,t_f]\times\d$ to $\R^\nx$. 
    Assume that $x_0 \in \M_0(t_0)$ is a given consistent initial value for \eqref{eq:semiDAE}
	at $t_0$. Moreover, Assumption~{\em\ref{ass:representation_y}} is valid for $\dx$.
	Additionally, for $\mathsf f$ the following conditions are fulfilled:
    \begin{enumerate}
		\item [\em 1)] $\mathsf f(t,\cdot\,,\cdot)$ is continuous on $\d$ f.a.a. $t \in [t_0,t_f]$;
		\item [\em 2)] for all $\mathsf z_\circ=(\mathsf x_\circ,\mathsf y_\circ) \in \d$ there are $\rho_\circ > 0$ and $\alpha_\mathsf f\in L^2(t_0,t_f)$ with
        \begin{align*}
            {\|f(t,\mathsf z)-f(t,\tilde{\mathsf z})\|}_2 \leq \alpha_\mathsf f(t)\,{\|\mathsf z-\tilde{\mathsf z}\|}_2\text{ f.a.a. }t\in(t_0,t_f],\text{ for all }\mathsf z,\tilde{\mathsf z}\in B_{\rho_\circ}(\mathsf z_\circ)\subset\d.
        \end{align*}
	\end{enumerate}
    Then there are $t_e\in(t_0,t_f]$ and a unique solution pair $(x,y): [t_0,t_e] \to \d$  (in the extended sense) to \eqref{eq:semiDAE} satisfying $x(t_0)=x_0$ and $y(t)=\mathsf y(t,x(t))$.
\end{theorem}

\begin{proof}
	We define the function $\mathsf F : [t_0,t_f] \times \dx \to \R^\nx$ by
    \begin{align}
        \label{Function_F}
        \mathsf F(t,\mathsf x):= \mathsf f(t,\mathsf x,\mathsf y(t,\mathsf x))\quad\text{f.a.a. }t\in[t_0,t_f]\text{ and for all }\mathsf x\in\dx.
    \end{align}
    It follows from condition 1) and Assumption~\ref{ass:representation_y}-4), -5) that the hypotheses 1) and 2) of Theorem~\ref{thm:caratheo_ex} are satisfied. Utilizing \eqref{Function_F}, condition 2) and Assumption~\ref{ass:representation_y}-6) we conclude that there exists $ \tilde\rho_\circ>0$ (small enough) with
    \begin{align*}
        &{\|\mathsf F(t,\mathsf x)-\mathsf F(t,\tilde{\mathsf x})\|}_2={\|f(t,\mathsf x,\mathsf y(t,\mathsf x))-f(t,\tilde{\mathsf x},\mathsf y(t,\tilde{\mathsf x}))\|}_2\\
        &\quad\le\alpha_\mathsf f(t)\,{\|(\mathsf x,\mathsf y(t,\mathsf x))-(\tilde{\mathsf x},\mathsf y(t,\tilde{\mathsf x}))\|}_2=\alpha_\mathsf f(t)\,\Big({\|\mathsf x-\tilde{\mathsf x}\|}_2^2+{\|\mathsf y(t,\mathsf x)-\mathsf y(t,\tilde{\mathsf x})\|}_2^2\Big)^{1/2}\\
        &\quad\le\alpha_\mathsf f(t)\,\Big({\|\mathsf x-\tilde{\mathsf x}\|}_2^2+\alpha_\mathsf y^2(t){\|\mathsf x-\tilde{\mathsf x}\|}_2^2\Big)^{1/2}\le\alpha_\mathsf f(t)\big(1+\alpha_\mathsf y(t)\big)\,{\|\mathsf x-\tilde{\mathsf x}\|}_2
    \end{align*}
    f.a.a. $t\in[t_0,t_f]$ and for all $\mathsf x,\tilde{\mathsf x}\in B_{\tilde\rho_\circ}(x_0)$. Due to $\alpha_\mathsf f,\alpha_\mathsf y\in L^2(t_0,t_f)$ we have $\alpha=\alpha_\mathsf f(1+\alpha_\mathsf y)\in L^1(t_0,t_f)$. Thus, condition 3) of Theorem~\ref{thm:caratheo_ex} holds. Finally, we infer from Assumption~\ref{ass:representation_y}-7) that hypothesis 4) of Theorem~\ref{thm:caratheo_ex} is satisfied as well. Therefore, there exist $t_e\in(t_0,t_f]$ and a unique solution 
	$x\in AC(\Ie;\mathbb R^\nx)$ to
	\begin{align}
        \label{DAE-x_eq}
		\dot x(t)&= \mathsf F(t,x(t))=\mathsf f(t,x(t),\mathsf y(t,x(t)))\text{ for } t \in \Ie,\quad x(t_0)=x_0
	\end{align}
	in the extended sense with $\Ie=[t_0,t_e]\subset[t_0,t_f]$. From $x\in AC(\Ie;\mathbb R^\nx)$ we infer that $x$ is continuous. Since $x(t_0)=x_0\in\dx$ holds and the set $\dx$ is open, we have -- possibly after choosing a smaller $t_e>t_0$ that $x(t)\in\dx$ for all $t\in\Ie$. Using \eqref{DAE-x_eq}, Assumption~\ref{ass:representation_y}-2) and -3) it follows from Theorem~\ref{thm:caratheo_ex} that the pair $(x,y)$ is a unique solution to \eqref{eq:semiDAE} in the extended sense satisfying $y(t)=\mathsf y(t,x(t))\in\dy$ on $\Ie$. Moreover, $y\in L^2(\Ie;\mathbb R^\ny)$ holds due to Assumption~\ref{ass:representation_y}-1).
\end{proof}
\begin{remark}\label{rem:bounded}
    Due to the continuity of the solution $x$ there exists a radius $\rho_\circ=\rho_\circ(x_0)>0$ such that -- 
    possibly after choosing a smaller $t_e>t_0$ -- we have $x(t) \in B_{\rho_\circ}(x_0)\subset\dx$ for $t \in \Ie$ hold. Thus,
    \begin{align*}
        {\|x\|}_{C(\Ie;\R^\nx)}&=\sup_{t\in\Ie}{\|x(t)\|}_2\le\sup_{t\in\Ie}\Big({\|x_0\|}_2+{\|x(t)-x_0\|}_2\Big)\le{\|x_0\|}_2+\rho_\circ=:c_\mathsf{apr},
    \end{align*}
    where the a-priori bound $c_\mathsf{apr}$ depends on the initial condition $x_0$.\hfill$\Diamond$ 
\end{remark}

\section{The optimization problem}
\label{sec:3}

In this section we focus on the optimization problem \eqref{eq:DAE_Orig} and prove existence of optimal solutions.

\subsection{The control-to-state operator}
\label{subsec:3.2}

We suppose that $x_0\in\xad \cap \M_0(t_0)$. Let us recall the DAE constraint of \eqref{eq:DAE_Orig} here: For $u\in\uad$ we have
\begin{align}
    \label{DAEstate}
    \begin{aligned}
        \dot{x}(t)&=f(t,x(t),y(t),u(t))&& \text{f.a.a. } t\in(t_0, t_f],\quad x(t_0)=x_0,\\
        0 &= g(t, x(t), y(t), u(t)) && \text{f.a.a. }t\in(t_0, t_f].
    \end{aligned}
\end{align}
Theorem~\ref{thm:existence_index1} ensures that \eqref{DAEstate} has a unique solution pair $(x,y)$ on $[t_0,t_e^u]$ for every $u\in\uad$, where $t_e^u\in(t_0,t_f]$ depends on $u$. Next we derive sufficient conditions that there is a $t_e\in(t_0,t_f]$ so that $t_e^u\ge t_e$  holds for any $u\in\uad$. For that purpose we suppose that the inhomogeinity $f$ has a specific structure, which is valid in our application.
\begin{assumption}
    \label{ass:represent_f}
    There exist open subsets $\dx\subset\R^\nx$, $\dy\subset\R^\ny$ and functions
    \begin{align*}
        f_1: [t_0,t_f]\times \dx \to \R^{\nx\times\ny},\quad f_2: [t_0,t_f]\times \dx \to \R^{\nx\times\nu},\quad f_3: [t_0,t_f] \times \dx \to \R^\nx
    \end{align*}
    satisfying 
    \begin{align*}
        f(t,\mathsf x,\mathsf y,\mathsf u) = f_1(t,\mathsf x)\mathsf y + f_2(t,\mathsf x)\mathsf u + f_3(t,\mathsf x)
    \end{align*}
    for $(t,\mathsf x,\mathsf y,\mathsf u) \in [t_0,t_f] \times \dx \times \dy \times\Uad$. The functions $f_1,$ $f_2$, $f_3$ and $f$ fulfill:
    \begin{enumerate}
        \item [1)] $f_1(\cdot\,,\mathsf x)$, $f_2(\cdot\,,\mathsf x)$ and $f_3(\cdot\,,\mathsf x)$ are measurable on $[t_0,t_f]$ for all $\mathsf x\in\dx$;
        \item [2)] $f_1(t,\cdot)$, $f_2(t,\cdot)$ and $f_3(t,\cdot)$ are continuous on $\dx$ f.a.a. $t \in [t_0,t_f]$;
        \item [3)] for all $\mathsf z_\circ=(\mathsf x_\circ,\mathsf y_\circ)\in\d:=\dx \times \dy$ there are $\rho_\circ>0$ and $\alpha_1,\alpha_2\in L^4(t_0,t_f)$ with 
        \begin{align*}
            {\|f(t,\mathsf z,\mathsf u)-f(t,\tilde{\mathsf z},\mathsf u)\|}_2&\leq \big(\alpha_1(t)\,{\|\mathsf u\|}_2+\alpha_2(t)\big){\|\mathsf z-\tilde{\mathsf z}\|}_2
        \end{align*}
        f.a.a. $t\in[t_0,t_f]$, for all $\mathsf z,\tilde{\mathsf z}\in B_{\rho_\circ}(\mathsf z_\circ)\subset\d$ and for all $\mathsf u\in\Uad$;
        \item [4)] for all $\mathsf x_\circ \in \dx$ there are $\rho_\circ>0$ and $\iota_j \in L^4(t_0,t_f)$ for $j = 1,2,3$ with 
        \begin{align*}
            {\|f_j(t,\mathsf x)\|}_F \leq \iota_j(t)\text{ for } j = 1,2,\quad{\|f_3(t,\mathsf x)\|}_2 \leq \iota_3(t)
        \end{align*}
        f.a.a. $t \in [t_0,t_f]$ and for all $\mathsf x \in B_{\rho_\circ}(\mathsf x_\circ)$.
    \end{enumerate}
\end{assumption}

To handle the algebraic function $g$, we need similar hypotheses as stated in Assumption \ref{ass:representation_y}, but now including the control function:
 
\begin{assumption}
    \label{ass:representation_yu}
    For \eqref{DAEstate} and the set $\uad$ of admissible controls there exists a function $\mathsf y\in L^2_{loc}((t_0,t_f)\times \dx\times\Uad;\R^\ny)$ with
    \begin{align*}
        g(t,\mathsf x,\mathsf y(t,\mathsf x,\mathsf u),\mathsf u)=0\quad\text{f.a.a. }t\in[t_0,t_f]\text{ and for all }(\mathsf x,\mathsf u)\in\dx\times\Uad.
    \end{align*}
    Additionally, there exist $\mathsf y_1:[t_0,t_f]\times\dx\to \R^{\ny\times\nu}$ and $\mathsf y_2:[t_0,t_f]\times\dx\to\R^\ny$ with 
    \begin{align*}
        \mathsf y(t,\mathsf x,\mathsf u)=\mathsf y_1(t,\mathsf x)\mathsf u+\mathsf y_2(t,\mathsf x)\quad\text{f.a.a. }t\in[t_0,t_f]\text{ and for all }(\mathsf x,\mathsf u)\in\dx\times\Uad.
    \end{align*}
    The representation fulfills the following conditions:
 	\begin{enumerate}
		\item [1)] $\mathsf y(t,\mathsf x,\mathsf u) \in \dy$ f.a.a. $t\in[t_0,t_f]$ and for all $(\mathsf x,\mathsf u)\in\dx\times[\ua,\ub]$;
		\item [2)] $\mathsf y_1(\cdot\,,\mathsf x)$ and  $\mathsf y_2(\cdot\,,\mathsf x)$ are measurable on $[t_0,t_f]$ for all $\mathsf x \in\dx$;
		\item [3)] $\mathsf y_1(t,\cdot)$ and $\mathsf y_2(t,\cdot)$ are continuous on $\dx$ f.a.a. $t \in [t_0,t_f]$;
		\item [4)] for all $\mathsf x_\circ \in \dx$ there are $\rho_\circ >0$ and $\gamma_1,\gamma_2\in L^4(t_0,t_f)$, with 
        \begin{align*}
            {\|\mathsf y_1(t,\mathsf x)-\mathsf y_1(t,\tilde{\mathsf x})\|}_F \leq \gamma_1(t)\,{\|\mathsf x-\tilde{\mathsf x}\|}_2,\\
            {\|\mathsf y_2(t,\mathsf x)-\mathsf y_2(t,\tilde{\mathsf x})\|}_2 \leq \gamma_2(t)\,{\|\mathsf x-\tilde{\mathsf x}\|}_2
        \end{align*}
        f.a.a. $t\in[t_0,t_f]$ and for all $\mathsf x,\tilde{\mathsf x}\in B_{\rho_\circ}(\mathsf x_\circ)\subset\dx$;
        \item [5)] for all $\mathsf x_\circ \in \xad$ there exist $\beta_1,\beta_2\in L^2(t_0,t_f)$ with 
        \begin{align*}
            {\|f(t,\mathsf x_\circ,\mathsf y(t,\mathsf x_\circ,\mathsf u),\mathsf u)\|}_2\leq \beta_1(t)\,{\|\mathsf u\|}_2+\beta_2(t)\quad\text{f.a.a. }t\in[t_0,t_f],\text{ for all }\mathsf u\in\Uad;
        \end{align*}
        \item [6)] for all $\mathsf x_\circ \in \dx$ there are $\rho_\circ$ and $\kappa_j \in L^4(t_0,t_f)$ with 
        \begin{align*}
            \|y_1(t,x)\|_F \leq \kappa_1(t) \mbox{ and } \|y_2(t,x)\|_2 \leq \kappa_2(t)
        \end{align*}
        f.a.a. $t \in [t_0,t_f]$ and for all $x \in B_{\rho_\circ}(\mathsf x_\circ)$.
	\end{enumerate}
\end{assumption}

\begin{remark}
    \label{rem:gnew}
    Depending on the choice for the control $u\in\uad$, the time interval, where \eqref{DAEstate} possesses a unique solution pair $(x,y)$, may change. Therefore, we need to derive bounds for terminal time that are independent of $u$.\hfill$\Diamond$
\end{remark}

The next result is proved in the appendix (see Section~\ref{Sec:A}).

\begin{lemma}
    \label{lem:bound_existencetime}
	Assume that $x_0\in\xad \cap \M_0(t_0)$, Assumptions~\ref{ass:represent_f} and \ref{ass:representation_yu} hold. Then there is a $t_e\in(t_0,t_f]$ such that for all $u \in \uad$ there exists a unique solution pair $(x,y)$ (in the extended sense) to \eqref{DAEstate} on $\Ie=[t_0,t_e]$. 
\end{lemma}

\begin{remark}\label{rem:a35}
    Let $x_0\in\x_0\cap\M_0(t_0)$. With Assumptions~\ref{ass:represent_f} and \ref{ass:representation_yu} holding it follows from Lemma~\ref{lem:bound_existencetime} that there exists a solution pair $(x,y)$ in the extended sense on the time interval $\Ie$ for every $u\in\uad$. Here, $\Ie$ does not depend on the chosen control function $u\in\uad$. With the arguments from the proof of Theorem \ref{thm:caratheo_ex} and from Remark~\ref{rem:bounded} we can suppose -- after possibly choosing a smaller $t_e>t_0$ -- that there exist a radius $\rho_\circ>0$ such that the solution $x$ fulfills $x(t) \in B_{\rho_\circ}(x_0)\in\dx$ for all $t \in\Ie$ and $u\in\uad$. Thus, we have
    \begin{align}
        \label{AprioriEst}
        \begin{aligned}
            {\|x\|}_{C(\Ie;\R^\nx)}\le c_\mathsf{apr}\text{ for every }u\in\uad,\quad
            {\|x(t)-x_0\|}_2\le\rho_\circ\text{ for all }t \in \Ie   
        \end{aligned}
    \end{align}
    with the a-priori constant $c_\mathsf{apr}=\rho_\circ+\|x_0\|_2$; cf. Remark~\ref{rem:goal}-1).\hfill$\Diamond$
\end{remark}

To ensure unique solvability on $[t_0,t_f]$ we make use of the next hypothesis.

\begin{assumption}
    \label{A-3}
    Suppose that $t_e=t_f$, i.e., $\Ie=[t_0,t_f]$ holds.
\end{assumption}

\begin{remark}
    \label{Rem:controlstate1}
    We define the spaces
    \begin{align*}
        \x:=AC([t_0,t_f];\R^\nx)\quad\text{and}\quad\y:=L^2([t_0,t_f];\R^\ny).
    \end{align*}
    Then, it follows from Lemma~\ref{lem:bound_existencetime} and Assumption~\ref{A-3} there exists a unique solution pair $(x,y)\in\x\times\y$ (in the extended sense) to \eqref{DAEstate} on $[t_0,t_f]$.\hfill$\Diamond$
\end{remark}



Let Assumptions \ref{ass:represent_f}, \ref{ass:representation_yu} and \ref{A-3} be valid. Due to Remark~\ref{Rem:controlstate1} the \emph{control-to-state operator}
\begin{align*}
    \mathcal S:\uad\to\x\times\y,\quad (x_u,y_u)=\mathcal S(u)\text{ solves \eqref{DAEstate} for given $u\in\uad$}
\end{align*}
is well-defined. For the proof of the next lemma we refer to the appendix (Section~\ref{Sec:B}).

\begin{lemma}
    \label{lem:controlstate2}
	Assume that Assumptions \ref{ass:represent_f}, \ref{ass:representation_yu} and \ref{A-3} hold. Let $x_0 \in \m_0(t_0)\cap\x_0$ be given and $\{u_k\}_{k \in \N} \subset \uad$ be a sequence in $\uad$. Consider the corresponding sequence $\{(x_k,y_k)\}_{k\in\N}$ with $(x_k,y_k)=\mathcal S(u_k)$. 
    Then there exists a subsequence 
	$\{x_{k_\ell}\}_{\ell \in \N}$ that converges to $\tx$ in $C([t_0,t_f];\R^\nx)$ with $(\tx,\ty) =\mathcal S(\tu)$ and $\ty=\mathsf y(\cdot\,,\tx(\cdot),\tu(\cdot))$. Furthermore, $\{x_{k_\ell}\}_{\ell \in \N}$ converges weakly to $\tx$ in $H^1(t_0,t_f;\R^\nx)$.
\end{lemma}

\subsection{Existence of optimal solutions}
\label{subsec:3.3}

Now we are able to prove the next result.

\begin{theorem}\label{thm:wellposed}
    Suppose that $x_0 \in \x_0 \cap \m_0(t_0)$ holds. Let
    Assumption~{\em\ref{ass:represent_f}}, {\em\ref{ass:representation_yu}} and {\em\ref{A-3}} hold. Assume 
    that $J : \x \times \y \times \uad \to \R$ is non-negative and weakly lower semi-continuous.
	Then there is a solution $(\tx,\ty,\tu)\in\x\times\y\times\uad$ to \eqref{eq:DAE_Orig}.
\end{theorem}

\begin{proof}
    Notice that the cost functional is non-negative. Let $\{(x_k,y_k,u_k)\}_{k\in\N}\subset\x\times\y\times\uad$ be a minimizing sequence for \eqref{eq:DAE_Orig}. Then, we have
    \begin{align*}
        \tJ=\inf\big\{J(x,y,u)\,\big|\,(x,y)=\mathcal S(u)\text{ and  }u\in\Uad\big\}=\lim_{k\to\infty}J(x_k,y_k,u_k).
    \end{align*}
	The set $\uad$ is a closed, convex, bounded subset of the Hilbert space $\u$. Thus $\uad$ is weakly compact, so that there exists a subsequence $\{u_{k_\ell}\}_{\ell\in\N}\subset\uad$ and an element $\tilde u\in\uad$ so that $u_{k_\ell}\rightharpoonup\tilde u$ in $\u$ as $\ell\to\infty$. We also have $(x_k,y_k)=\mathcal S(u_k)$ and $y_k(t)=\mathsf y(t,x_k(t),u_k(t))$ f.a.a. $t\in[t_0,t_f]$ for all $k\in\N$. It follows from Lemma~\ref{lem:controlstate2} that we can assume -- possibly after taking a further subsequence -- that $\{x_{k_\ell}\}_{\ell\in\N}\subset\x$ converges uniformly to the element $\tx$ with $(\tx,\ty)=\mathcal S(\tu)\in\x$ and $\ty=\mathsf y(\cdot\,,\tx(\cdot),\tu(\cdot))$. We have
    \begin{align*}
        \mathsf y_j(t,x_{k_\ell}(t)) \to \mathsf y_j(t,\tx(t))&\quad\text{f.a.a. }t \in [t_0,t_f]\text{ and }j=1,2,
    \end{align*}
    due to Assumption~\ref{ass:representation_yu}-3). 
    With Assumption~\ref{ass:representation_yu}-6) we can apply the dominated convergence theorem \cite[p.~321]{Rud76} and conclude that the sequences $\{\mathsf y_1(\cdot,x_{k_\ell}(\cdot))\}_{\ell\in\mathbb N}$ and $\{\mathsf y_2(\cdot,x_{k_\ell}(\cdot))\}_{\ell\in\mathbb N}$
    converge in $L^2(t_0,t_f;\R^{\ny\times\nu})$ and $L^2(t_0,t_f;\R^\ny)$, respectively.
    Then with $u_{k_\ell} \rightharpoonup\tu$ in $\u $ for $\ell\to\infty$ we can apply \cite[Theorem 5.12-4 (c)]{Ciarlet13} and conclude that 
    $y_{k_\ell} \rightharpoonup\ty$ in $\y $ for $\ell\to\infty$.
    Again applying Lemma~\ref{lem:controlstate2} we derive that -- possibly after taking a further subsequence -- $\{x_{k_\ell}\}_{\ell\in\mathbb N}$ converges weakly to $\tx$ in $H^1(t_0,t_f;\R^\nx)\hookrightarrow \x$. Thus, we have $(x_{k_\ell},y_{k_\ell},u_{k_\ell})\rightharpoonup(\tx,\ty,\tu)$ for $\ell\to\infty$ in $\x\times\y\times\u$. Hence, the weakly lower semi-continuity of $J$ guarantees
	\begin{align*}
		\tJ = \lim\limits_{\ell \to \infty} J(x_{k_\ell},y_{k_\ell},u_{k_\ell}) = \liminf\limits_{\ell \to \infty}  J(x_{k_\ell},y_{k_\ell},u_{k_\ell}) 
            \geq J(\tx,\ty,\tu) \geq \tJ,
	\end{align*}
	which implies that $(\tx,\ty,\tu)$ is a minimizer for \eqref{eq:DAE_Orig}.
\end{proof}

\section{Network Optimization}
\label{sec:4}

In this section we introduce a district heating network model which can be expressed like \eqref{eq:semiDAE}. For more details we refer the reader to \cite{JaekleReichle} and, e.g., to \cite{Borsche19, Koecher00, Krug19}. 

In a district heating network the hot water is distributed to households by a system of pipes. A network with a not necessarily identical structure transports the colder water back to the power plant. In each pipe we describe the flow by the following three equations:
\begin{subequations}
    \label{eq:endeq}
    \begin{align}
        \label{eq:transport}
        \partial_t\rho + \partial_x (\rho v)&= 0,\\
        \partial_t (\rho v) + \partial_x\left(\rho v^2 + p\right)&= -\frac{\lambda}{2d}v|\rho v| -G(\partial_x h)\rho,\\
        \label{eq:transport_temperature}
        \partial_t(c_p\rho T) + \partial_x(\rho v c_pT)&= -\frac{4k}{d}(T-T_{ext}).
    \end{align}
\end{subequations}
These three equations are the conservation of mass, the balance of 
momentum and the balance of energy. The term $G(\partial_xh)$ takes the 
vertical displacement $h$ into account, where $G$ is the gravitational acceleration. In a pipe with length $L$ and height difference $\Delta h$ it holds 
$G(\partial_x h) = G \Delta h /L$ with $G \approx 9.80665 ms^{-2}$. All variables and parameters are summarized in Table~\ref{tab:Var}.
\begin{table}[!t]
\caption{Variables (top) and parameters (bottom) of the district heating network model.}
\label{tab:Var}       
\begin{center}
{\small\begin{tabular}{p{2.3cm}p{5.8cm}p{2.5cm}p{1.8cm}}
\hline\hline \noalign{\smallskip}
Symbol & Explanation & Unit & Example    \\
\hline \noalign{\smallskip} \noalign{\smallskip}
$v(t)$ & Flow velocity  & m s$^{-1}$ & - \\
$p(t,0)$, $p(t,L)$ & Pressure at the ends of a pipe & Pa$=$kg m$^{-1}$s$^{-2}$& $5\cdot 10^5, 2\cdot 10^5$\\
$T(t,x),T_{i,j}(t)$  & Water temperature in a pipe  & C& - \\
$P_p(t)$ & Pumping power & W & - \\
$P_w(t)$ & Waste incineration power & W & - \\
$P_g(t)$ & Gas combustion power & W & - \\
\noalign{\smallskip}\hline\noalign{\smallskip}
$\rho$ & Density of the water & kg m$^{-3}$& $960$ \\
$t$ & Time coordinate; $t \in [t_0,t_f]$  & s& -\\
$x$ & Spatial coordinate in a pipe & m &-\\
$L$ & Length of a pipe & m &$300$\\
$\Delta h$ & Height difference in a pipe & m & - \\
$d$ & Diameter of a pipe & m &$0.1$\\
$A$ & Cross-sectional area of a pipe; & m$^{2}$ & $7.85\cdot 10^{-3}$\\
$\lambda$ & Friction factor of a pipe & $1$&$0.00251$ \\
$Q_k(t)$ & Power consumption & W & $1\cdot 10^6$\\
$k_{rough}$ & Roghness of the inner wall & m &$0.00026$\\
$k$ & Heat transfer coefficient & W m $^{-2}$ C $^{-1}$ &$0.1$\\
$\bar T_\mathsf{out}$ & Water temperature & C& $60$\\
$T_\mathsf{ext}$ & Surrounding temperature & C& $20$\\
$c_p$ & Heat capacity of water & J kg$^{-1}$C$^{-1}$ &$4160$\\
$G(\partial_x h)$ & Gravitational acceleration  & m s$^{-2}$ &$9.80665$\\
$\omega_w, \omega_g,\omega_p$ & Cost coefficient & \euro{}/W & $4,5$ \\
$\alpha$, $\beta$, $\gamma$ & Cost coefficients & - & - \\
$N$; $n_q$; $\bar N$ & Number of pipes; consumers; nodes & - & - \\
\noalign{\smallskip}\hline\noalign{\smallskip}
\end{tabular}}
\end{center}
\end{table}
We assume that our network has $N$ pipes and $I = \{1,...,N\}$ describes the index set for the pipes. For the index set $I$ we use two different representations. In the first one we distinguish whether $i \in I$ is in the forward ($i\in I_f$) or in the backward flow ($i\in I_b$). Therefore, we have $I = I_f \cup I_b$. In the second one we can divide the index set $I$ into three sets, depending if pipe $i \in I$ is a supply ($i\in I_S$), demand  ($i\in I_D$) or interior pipe ($i\in I_I$). 
Therefore, we have the representation   $I = I_S \cup I_D \cup I_I$.
The number of consumsers is $n_q$ and $\bar N = n_s + n_d + n_{junc}$
the number of nodes, while $n_s$ is the number of supply nodes, $n_d$ the number of demand nodes and $n_{junc}$ the number of interior nodes.

Next we simplify the system, whereby, for example. the density $\rho$ 
is assumed to be constant. At the beginning, the density $\rho = \rho(t,x)$, 
the velocity $v = v(t,x)$ and the temperature $T = T(t,x)$ describe time and location-dependent functions. For better readability we often use $\rho$, $v$ or $T$ instead of $\rho(t,x)$, $v(t,x)$ or $T(t,x)$. 

To model heating network we need coupling conditions at the junctions additional to the flow equations on the edges. Therefore, we consider the following coupling conditions in every interior node $j$:
\begin{subequations}
\label{couplingCond}
\begin{align}
\label{conserveMass}
\sum_{i \in \sigma_j} A_i \rho_i(t,L_i) v_i(t,L_i) 
    &= \sum_{i\in \Sigma_j}A_i \rho_i(t,0) v_i(t,0),  \\
\label{conserveEnergy}
\sum_{i \in \sigma_j} c_p A_i \rho_i(t,L_i) v_i(t,L_i) T_i(t,L_i) 
    &= \sum_{i\in \Sigma_j}c_p A_i \rho_i(t,0) v_i(t,0) T_i(t,0), \\
\label{contPressure}
p_i(t,L_i) &= p_l(t,0), \quad \text{for all } i \in \sigma_j, \, l\in \Sigma_j, \\
\label{perfectMix}
T_i(t,0) &= T_l(t,0), \quad \text{for } i,l \in \Sigma_j,\, i \neq l,
\end{align}
\end{subequations}
where
\begin{itemize}
    \item $\sigma_j$ is the set of all incoming pipes at node node $j$ and
    \item $\Sigma_j$ is the set of all leaving pipes at node $j$.
\end{itemize}
Equation \eqref{conserveMass} states the conservation of mass and \eqref{conserveEnergy} the conservation of energy. The continuity of the pressure -- described by \eqref{contPressure} -- is a widely used condition, see, e.g., \cite{Banda06, Colombo08, Domschke15}. Additionally we assume a perfect mixing of flows at the junction, which means that we assume the same temperature in all outgoing pipes; cf. \eqref{perfectMix}.

An additional important component in a district heating network are the consumers which connect the supply network with the return one. Each consumer needs a certain amount of thermal power $Q_k(t)$ at time $t$ with $k = 1,...,n_q$. We use the index ``$\mathsf{in}$'' for the incoming pipe and ``$\mathsf{out}$'' for the outgoing pipe before and after a consumer. The outgoing temperature $\overline{T}_\mathsf{out}$ is assumed to be a fixed value and no mass is lost. This leads to the following consumer equations:
\begin{subequations}
\label{connectNetwork}
\begin{align}
\label{noMassLost}
\rho_{in}(t,L_\mathsf{in}) v_\mathsf{in}(t,L_\mathsf{in}) &= \rho_\mathsf{out}(t,0) v_\mathsf{out}(t,0),\\
\label{consumerDemand}
Q_k(t) &= c_p A \rho_\mathsf{in}(t,L_\mathsf{in}) v_\mathsf{in}(t,L_\mathsf{in})(T_\mathsf{in}(t,L_\mathsf{in}) - \overline{T}_\mathsf{out}),
\end{align}
\end{subequations}
where $T_\mathsf{in}$ is the temperature of the flow arriving at the household.

Next we introduce the power plant, where we add the control function $u$. 
Let us assume that the power plant is located between the first and the last pipe 
in the network and we always have one power plant in the system. We assume (cf. \cite{Krug19}) that the power to run the pumps can be controlled in such a way that a desired pressure increase in the depot of the district heating network can be realized. This control is denoted by $P_p$. Moreover a temperature gain is obtained by thermal power production in the depot. The functions $P_w$ and $P_g$ describe the thermal power produced by waste 
incineration and gas combustion. We have
\begin{subequations}
\label{eq:controlFun}
  \begin{align}
    \label{eq:controlFun1}
      P_p(t) &= A_1 v_1(t,0) \left(p_1(t,0) - p_N(t,L_N) \right),\\
    \label{eq:controlFun2}
    P_w(t) + P_g(t) &= A_1 \rho_1(t,0) v_1(t,0)c_p \left( T_1(t,0) - T_N(t,L_N) \right).
  \end{align}
\end{subequations}
Note that the second equation is similar to the consumer equation \eqref{consumerDemand}. The water arrives with a given pressure $p_{d}$ at the depot, the so-called pressure stagnation which we need in this model formulation to achieve a unique pressure distribution in the network: 
\begin{equation}\label{eq:pressure_stagnation}
    p_N(t,L_N) = p_{d}(t).
\end{equation}

In order to reduce the complexity of the system, we assume some simplifications which, 
in our case, produce negligible errors. As the water in the pipes is almost incompressible and the temperature difference is small 
we suppse $\partial_t \rho = 0$ and $\partial_x v = 0$.

Another assumption is that the network 
has a tree structure and all consumers have a minimum consumption of $1\,W$. In this setting we can assume that the sign of the velocity is always positive and simplify $|\rho v|v$ to $\rho v^2$. We neglect the term $\partial_t v$ as described in \cite{JaekleReichle} since this in the end leads to a system with differentiation index $di=1$; cf. also \cite{hauschild2020port}. In \cite{JaekleReichle} we have discussed in more detail that this assumption is suitable in a district heating network.

In the last step we apply a spatial discretization by the methods of lines for the partial differential equation \eqref{eq:transport_temperature} and get 
\begin{equation}
    \label{eq:temperature_disc1}
    \dot T_{i,j}(t) + \frac{v_i(t)}{\Delta x}(T_{i,j}(t)-T_{i,j-1}(t)) 
    = -\frac{4k}{c_pd\rho}(T_{i,j}(t) - T_\mathsf{ext})
\end{equation}
for $j = 2,...,n_i$ and $i = 1,...,N$.
Summarizing, for a system with $N$ pipes, $n_q$ consumers and $\Bar{N} = n_s + n_d + n_{junc}$ nodes we receive the following 
DAE system: For every pipe $i \in I$ and f.a.a. $t\in[t_0, t_f]$ we have:
\begin{subequations}
    \label{eq:index1system_pipes}
    \begin{align}
        \dot T_{i,j}(t)
            &= -\frac{v_i(t)}{\Delta x}(T_{i,j}(t)-T_{i,j-1}(t))-\frac{4k}{c_pd\rho}(T_{i,j}(t)-T_\mathsf{ext})\quad\text{for } j = 2,...,n_i,
            \label{eq:index1system_temp}\\
            0 &= p_i(t,L_i)-p_i(t,0)+\frac{\lambda\rho L}{2d}\,v_i(t)^2 + LG(\partial_x h)\rho.
            \label{eq:index1system_velo}
    \end{align}
For every consumer $k$ we have:
    \begin{align}
        0 
            &= v_\mathsf{in}(t)-v_\mathsf{out}(t) 
                \label{eq:index1system_cons1},\\
         0 
            &= c_pA_\mathsf{in}\rho v_\mathsf{in}(t) (T_{\mathsf{in},n_\mathsf{in}}(t)-\overline{T}_\mathsf{out})-Q_k(t)
                \label{eq:index1system_cons2},\\
        0
            &= T_{\mathsf{out},1}(t) - T_\mathsf{out}(t).
                \label{eq:index1system_cons3}
    \end{align}
For every interior node $j \in I_I$ we have:
    \begin{align}
        0 
            &= \sum_{i \in \sigma_j} A_i v_i(t) - \sum_{i\in \Sigma_j}A_iv_i(t)
                \label{eq:index1system_consmass}, \\
        0
            &= \sum_{i \in \sigma_j} c_p A_i v_i(t) T_{i,n_i}(t)-\sum_{i\in \Sigma_j}c_p A_i v_i(t) T_{i,1}(t)
                \label{eq:index1system_consenerg},\\
       0 
            &=  p_i(t,L_i) - p_l(t,0)
                \mbox{ for all } i \in \sigma_j, \, l\in \Sigma_j 
                \label{eq:index1system_conti},\\
        0
            &= T_{i,1}(t)  - T_{l,1}(t) \mbox{ for } i,l \in \Sigma_j,\, i \neq l.
                \label{eq:index1system_perfectmix}
    \end{align}
For the control and the pressure in the end we have:
    \begin{align}
        0 
            &= P_p(t)-A_1v_1(t)(p_1(t,0)-p_N(t,L_N)),
            \label{eq:index1system_control1}\\
        0
            &= P_w(t)+P_g(t) - A_1\rho v_1(t)c_p(T_{1,1}(t)-T_{N,n_N}(t)),
            \label{eq:index1system_control2}\\
        0 
            &= p_N(t,L_N) - p_d(t).
            \label{eq:index1system_pressure}
    \end{align}
\end{subequations}
We build the differential variable $x$, the algebraic variable $y$ 
and the control variable $u$ in the following way: 
\begin{align*}
   x(t) &= \begin{pmatrix}
            T_{1,2}(t), \hdots
            ,T_{1,n_1}(t),
            \hdots
            ,T_{N,2}(t),
            \hdots
            ,T_{N,n_N}(t)
         \end{pmatrix}^\top, \\
    y(t) &=  \begin{pmatrix}
             v_1(t),
             \hdots
             ,v_N(t),\,
             T_{1,1}(t),
             \hdots
             ,T_{1,N}(t),\,
             p_{1}(t,0),
             \hdots 
             ,p_{N}(t,L_N)
         \end{pmatrix}^\top, \\ 
    u(t) &= \begin{pmatrix}
             P_p(t),\,
             P_w(t),\,
             P_g(t)
         \end{pmatrix}^\top
\end{align*}
with which we can describe \eqref{eq:index1system_pipes} as a DAE with differentiation index $di=1$. We define $f: [t_0,t_f] \times \dx \times \dy \to \R^\nx$ with $\nx = \sum_{i = 1}^N n_i-1$, $\ny=4N$ and $\nu=3$. Let
\begin{itemize}
    \item $g_1 : [t_0,t_f] \times \dx \times \dy \to \R^N$ through 
        \eqref{eq:index1system_cons1}, \eqref{eq:index1system_cons2} and \eqref{eq:index1system_consmass},
    \item $g_2 : [t_0,t_f] \times \dx \times \dy \times \uad \to \R^{3N}$ through \eqref{eq:index1system_velo}, \eqref{eq:index1system_cons3}, 
        \eqref{eq:index1system_consenerg}-\eqref{eq:index1system_pressure}. 
\end{itemize}
Summarizing, we get
\begin{subequations}\label{eq:our_dae}
    \begin{align}
        \Dot{x}(t) &= f(t,x(t),y(t)) &&\text{f.a.a. }t\in(t_0,t_f],\\
        0 &= g_1(t,x(t),y(t)) &&\text{f.a.a. }t\in(t_0,t_f],\\
        0 &= g_2(t,x(t),y(t),u(t))&&\text{f.a.a. }t\in(t_0,t_f]
    \end{align}
\end{subequations}
with
\begin{align*}
    \dx &:= \big\{\mathsf x\in\R^\nx_{>0}:\mathsf x_i>\overline{T}_\mathsf{out}\text{ for }i \in I_f \big\},\\
    \dy &:= \big\{\mathsf y \in \R^\ny : \mathsf y_i > 0 \mbox{ for } i = 1,...,2N\big\},\\
    \Uad &:= \big\{ u \in \R^3 : 0 \leq u_i \leq \ub\big\}\qquad\text{with }\ub \in \R_{> 0}.
\end{align*}

\begin{lemma}\label{lem:representation}
	Our given system \eqref{eq:our_dae} fulfils Assumptions 
    \ref{ass:represent_f} and \ref{ass:representation_yu} with $\dx, \, \dy$ and 
    $\uad$ defined as before. 
    Therefore, we have a representation of the algebraic variable $y$ depending 
    on $t$, $x$ and $u$.
\end{lemma}
\begin{proof}
    If $i \in I_D\backslash \{N\}$ we can 
    use \eqref{eq:index1system_cons2} and get 
    \begin{align}\label{eq:lemrep_vi}
        v_i(t) = \frac{Q_k(t)}{(T_{i,n_i}(t)-\overline{T}_\mathsf{out})c_pA_i\rho}.
    \end{align}
    The definition of $\dx$ implies, that $v_i$ is well-defined and if $Q_k(t) > 0$ holds, 
    we have $v_i(t) > 0$. \\
    If $i \in I_S\backslash \{1\}$ we can 
    apply \eqref{eq:index1system_cons1} and use the corresponding demand pipe $j_i \in I_D$ with \eqref{eq:lemrep_vi} to get 
    \begin{align}\label{eq:lemrep_vi2}
        v_i(t) = v_{j_i}(t) = \frac{Q_k(t)}{(T_{j_i,n_{j_i}}(t)-\overline{T}_\mathsf{out})c_pA_{j_i}\rho}
    \end{align}
    and again, if $Q_k(t) > 0$ holds, 
    we have $v_i(t) > 0$. \\
    If pipe $k$ is in the forward flow part of the network and an interior pipe or 
    the first pipe 
    we know (due to the tree structure) that the corresponding node $j$, which is located 
    at the end of pipe $k$, only has one incoming pipe. 
    Therefore \eqref{eq:index1system_consmass} is simplified to 
    \begin{align*}
        v_k(t) = \sum\limits_{i\in \Sigma_j} \frac{A_iv_i(t)}{A_k}
    \end{align*}
    and it holds $v_k(t) > 0$. \\
    If pipe $k$ is in the backward flow part of the network and an interior pipe or the last pipe, we get with the same argumentation a representation 
    \begin{align*}
        v_k(t) = \sum\limits_{i\in \sigma_j} \frac{A_iv_i(t)}{A_k}
    \end{align*}
    with $j$ is the node located in the beginning of pipe $k$. It holds $v_k(t) > 0$. 
    Iteratively with \eqref{eq:lemrep_vi}, \eqref{eq:lemrep_vi2} and the tree-structure we 
    have a representation of $v(t) = (v_1(t),...,v_N(t))^\top$ only depending on $x$, 
    for $x\in\dx$ and with $v(t) > 0$.\\
    For the temperature $T_{1,1}(t)$ we use \eqref{eq:index1system_control1} to get 
    \begin{align}\label{eq:lemrep_t11}
        T_{1,1}(t) = \frac{P_w(t)+P_g(t)}{A_1\rho v_1(t)c_p}+T_{N,n_N}(t)
    \end{align}
    which is well-defined for $x \in \dx$, $u \in \uad$ and it holds $T_{1,1}(t) > 0$. \\
    For pipe $i$ with $i \in I_S\backslash\{1\}$ we can use \eqref{eq:index1system_cons3}
    to get $T_{i,1}(t) = \overline{T}_\mathsf{out} > 0.$\\
    For pipe $i$ with $i \in I_D \cup I_I$ we use \eqref{eq:index1system_consenerg} and 
    \eqref{eq:index1system_perfectmix} to get 
    \begin{align}\label{eq:lemrep_tij}
        T_{i,1}(t) = \frac{\sum\limits_{i \in \sigma_j}c_pA_iv_i(t)T_{i,n_i}(t)}{\sum\limits_{i \in \Sigma_j}c_pA_iv_i(t)}
    \end{align}
    where $j$ is the corresponding node in the beginning of pipe $i$. If we use the 
    representation of $v_i(t)$, $x \in \dx$ and $u \in \uad$ we get a representation 
    of $T_{i,1}(t)$ for $i \in I$ independent of $y$ which fulfil $T_{i,1}(t) > 0$.\\
    For pipe $N$ and $1$ we have with \eqref{eq:index1system_pressure},
    \eqref{eq:index1system_control1} and \eqref{eq:index1system_velo} 
    the representations 
    \begin{align*}
        &p_N(t,L_N) = p_d(t), &p_N(t,0) = p_d(t) + \frac{\lambda}{2d}\rho L_N v_N(t)^2 + L g(\partial_x h)\rho, \\
        &p_1(t,0) = \frac{P_p(t)}{A_1v_1(t)}+p_d(t),
        &p_1(t,L_1) = p_1(t,0) - \frac{\lambda}{2d}\rho L_1 v_1(t)^2 - L G(\partial_x h)\rho.
    \end{align*}
    For pipe $i \in I_f\backslash\{1\}$ we use \eqref{eq:index1system_conti}, 
    \eqref{eq:index1system_velo} and the corresponding $j_i \in I_f$ with 
    \begin{align*}
        p_i(t,0) = p_{j_i}(t,L_{j_i}),\quad
        p_i(t,L_i) = p_{j_i}(t,L_{j_i}) 
            - \frac{\lambda}{2d}\rho L_i v_i(t)^2 - L G(\partial_x h)\rho.
    \end{align*}
    Then, iteratively the representations of $p_i(t,0)$
    and $p_i(t,L_i)$
    for $i \in I_f$, due to the tree-structure and the representation of $v(t)$ independent of $y$, 
    are independent of $y$.
    For pipe $i \in I_b\backslash\{N\}$ we use \eqref{eq:index1system_conti}, 
    \eqref{eq:index1system_velo} and the corresponding $j_i \in I_b$ with 
    \begin{align*}
        p_i(t,L_i) &= p_{j_i}(t,0)\\
        p_i(t,L_i) &= p_{j_i}(t)+\frac{\lambda}{2d}\rho L_i v_i(t)^2 + L G(\partial_x h)\rho.
    \end{align*}
    As before we have a representation independent of $y$.
    In summary, we receive a representation $\mathsf y(t,\mathsf x,\mathsf u) = \mathsf y_1(t,\mathsf x)\mathsf u+\mathsf y_2(t,\mathsf x)$ 
    which fulfills $g(t,\mathsf x,\mathsf y(t,\mathsf x,\mathsf u),\mathsf u) = 0$ for $(t,\mathsf x,\mathsf u) \in [t_0,t_f]\times\dx\times\uad$.
    In the next step we take a look at the function $f$. 
    Let $i \in \{1,...,N\}$ and $j \in \{1,...,n_i\}$. \\
    If $i = 1,\,j = 2$ we use \eqref{eq:lemrep_t11} and modify \eqref{eq:index1system_temp}
    to get 
    \begin{align*}
        \dot T_{1,2}(t) = -\frac{T_{1,2}(t)-T_{N,n_N}(t)}{\Delta x} v_1(t)
        +\frac{P_w(t)+P_g(t)}{\Delta x A_1\rho c_p}-\frac{4k}{c_p d\rho}(T_{1,2}(t)-T_{ext}).
    \end{align*}
    If $i > 1$ and $j = 2$ we use the representation \eqref{eq:lemrep_tij} of $T_{i,1}(t)$
    and the representation of $v(t)$ to modify \eqref{eq:index1system_temp}.
    Therefore we have
    \begin{align*}
        f(t,\mathsf x,\mathsf y,\mathsf u) = f_1(t,\mathsf x)\mathsf y + f_2(t,\mathsf x)\mathsf u + f_3(t,\mathsf x).
    \end{align*}
    For $\mathsf x\in\dx$ and $\mathsf u \in \uad$ we have $\mathsf y(t,\mathsf x,\mathsf u) \in \dy$. The functions 
    $\mathsf y_1(\cdot\,,\mathsf x)$, $\mathsf y_2(\cdot\,,\mathsf x)$, $f_1(\cdot\,,\mathsf x)$, $f_2(\cdot\,,\mathsf x)$ and $f_3(\cdot\,,\mathsf x)$
    are continuous and therefore measurable for each fixed $\mathsf x \in \dx$. 
    The functions $f_j$ and $y_i$ for $j = 1,2,3$ and $i = 1,2$ are once continuously differentiable. Therefore for every $x_0 \in \dx$ there exist $\rho_\circ > 0$ and functions $\iota_i,\,\kappa_j \in L^4(t_0,t_f)$, $i = 1,2,3$ and $j = 1,2$, satisfying
    \begin{align*}
        {\|f_i(t,\mathsf x)\|}_F\leq \iota_i(t)\quad\text{and}\quad{\|\mathsf y_1(t,\mathsf x)\|}_2\leq \kappa_1(t)\mbox{ for } i = 1,2\\
        {\|f_3(t,\mathsf x)\|}_F\leq \iota_3(t)\quad\text{and}\quad{\|\mathsf y_2(t,\mathsf x)\|}_2\leq \kappa_2(t)
    \end{align*}
    for all $t \in [t_0,t_f]$ and $\mathsf x \in B_{\rho_\circ}(x_0)$.
    In addition, the 
    functions $y_1(t,\cdot)$, $y_2(t,\cdot)$, $f_1(t,\cdot)$, $f_2(t,\cdot)$ and $f_3(t,\cdot)$
    are continuous. 
    The locally Lipschitz condition for $\mathsf y_1$, $\mathsf y_2$, $f_1$, $f_2$ and $f_3$ in $\mathsf x$ follow, 
    because the functions are once continuously differentiable and the interval $[t_0,t_f]$
    is bounded. 
    For the last condition in Assumption \ref{ass:representation_yu} we fix $x_0 \in \xad \cap \m_0(t_0)$ and that $f$ and $y$
    are continuous and choose $\beta_1(t)$ depending on $Q_k(t)$ with $k = 1,...,n_q$.
\end{proof}

The different variables of the network used in the model are subject to the 
following state and control constraints:
\begin{align*}
x(t) \in [x_\mathsf a,x_\mathsf b],\quad y(t) \in [y_\mathsf a,y_\mathsf b],\quad u(t) \in [u_\mathsf a, u_\mathsf b].
\end{align*}
The objective function to minimize is given by the linear cost
\begin{align*}
\int_{t_0}^{t_f} \left( \omega_w(t) P_w(t) + \omega_g(t) P_g(t) + \omega_p(t) P_p(t)  \right)\,\mathrm dt 
=: \int_{t_0}^{t_f} \underbrace{\sum_{i=1}^3 \omega_i(t) u_i(t)}_{=:f_0(u(t))}\,\mathrm dt=:J(u), 
\end{align*}
where $\omega_w$, $\omega_w$ and $\omega_p$ are positive cost coefficients of the 
waste incineration, the gas combustion and the pumping power, respectively. 
All in all, this leads to an OCP where the cost only depends on the control $u$. 
Without an initial value $x_0$ and boundary constraints on the control $u$, 
the state $x$ and the algebraic variables $y$ which is hard to solve and not in the 
form of Problem \eqref{eq:DAE_Orig}. Therefore we do some simplifications. 
First of all we use a regularization approach and change the cost to
\begin{align}\label{eq:our_cost}
\begin{aligned}
    &\int_{t_0}^{t_f} \sum_{i=1}^3 \omega_i(t) u_i(t) + 
    \frac{\alpha}{2} \Vert x(t) - x_d(t) \Vert^2_2 + 
    \frac{\beta}{2} \Vert y(t) - y_d(t) \Vert^2_2 \\
    &\hspace{8mm}+ \frac{\gamma}{2} \Vert u(t)-u_d(t) \Vert^2_2\,\mathrm dt
\end{aligned}
\end{align}
with some regularization parameters $\alpha, \beta, \gamma$ and functions $x_d$, $y_d$ and $u_d$. 
To get a suitable initial guess for the optimization we use an instantaneous control approach. 
Basically, the idea is to discretize the time interval and solve many stationary problems. Then we stack the solution to one control together which is usually far away from an optimal
solution but a good initial guess. Instantaneous control is closely related to receeding horizon control (RHC) or model predictive control (MPC) with finite time horizon \cite{instant1, instant2, instant3, instant4}. Its difference to RHC and MPC can best be explained as follows (cf. Bewley et al. in \cite{instant5}). (RHC) and (MPC) computes the next optimal controls and apply the first one to steer the system to the next time instance. Instantaneous control computes an approximation to the next optimal move by taking only into account the next time instance (instead of a larger finite horizon) to get an optimal control for the next time instance. It cannot be expected that the instantaneous control strategy stabilize arbitrary dynamical systems but it is capable -- as we will see in our numerical test - to provide a good starting value for the computation of an optimal control on $[t_0, t_f]$.
A pseudocode for an instantaneous control approach is given in Algorithm~\ref{alg:Instant}.
\begin{algorithm}
\caption{(Instantaneous Control Approach)}
\label{alg:Instant}
\begin{algorithmic}
\Require Discretize time Interval $\mathcal{T}$ with discrete timesteps 
$t_0=t_1 < t_2 < \ldots < t_{n_t} = t_f$ and $\Delta t_i = t_i - t_{i-1}$
\For{$i=1,2,\ldots,n_t$}
\State Compute the optimal solution $u_i^*$ tothe stationary problem (e.g., by a trust region method)\\
	\vspace{0.1cm}
		\begin{align*}
			\min_{u_i \in [u_\mathsf a,u_\mathsf b]} f_0(u_i)\quad\text{subject to}\quad\left\{
   		    \begin{aligned}
			     0 &=x_i - x_{i-1}- \Delta t_if(t_i, x_i, y_i, u_i) \\
			     0 &= g(t_i, x_i, y_i, u_i)
		    \end{aligned}
            \right.
        \end{align*}
\EndFor
\State Concatenate solution $u^* = (u_1^*, \ldots, u_N^*)$
\end{algorithmic}
\end{algorithm}

\begin{remark}\label{rem:desired}
The solution of the first time step is now used as a consistent initial value $x_0$ and the stacked 
solution of the state and algebraic variable is used for the the desired variables $x_d$ 
and $y_d$. For physical reasons and as a final simplification we assume that the 
restriction to the input $\ua \leq u \leq \ub$ and the regularization in the cost are sufficient to ensure the state constraints.\hfill$\Diamond$     
\end{remark}
\begin{remark}
    Note, that the optimal control problem \eqref{eq:DAE_Orig} based on \eqref{eq:our_dae} and \eqref{eq:our_cost} fulfills Assumptions~ \ref{ass:represent_f} and \ref{ass:representation_yu} from Theorem \ref{thm:wellposed}.\hfill$\Diamond$
\end{remark}

\section{Numerical Experiments}
\label{sec:5}
In this section we present numercial experiments for different networks:
\begin{enumerate}
    \item We consider a small network with a known analytic solution and compare the error of the analytic and the numerical optimal solution for different fine discretizations. It is shown that our method yields sufficiently good results in acceptable computational time. 
    \item In the second example we deal with a larger network and compare the numerical results while varying the different parameters in the cost function.
    \item Finally we compute an optimal solution for a real network with realistic data from our project partners.
\end{enumerate}
Direct discretization methods are often preferred for solving optimal control problems with real applications. Often they are able to return results even for difficult problems even with control and state constraints; cf., e.g., \cite{betts2001practical}. But from a mathematical point of view it would be better to investigate the optimal control problem subject to a DAE as an infinite dimensional optimal control problem and derive the optimality system. This approach leads to a nonlinear multi-point boundary value problem which usually yields the most accurate solutions. However, it suffers from the drawback that it requires a good initial guess and carefully chosen discretization parameters in order to converge and even then convergence is difficult to guarantee in applications. In, e.g., \cite{martens2019convergence, martens2021convergence} the authors show under certain assumptions that the discretized (locally) optimal control converges to the infinite dimensional one. For this reason we use a \emph{Discretize-before-Optimize} (DBO) approach for all numerical tests, where we discretize the DAE with an implicit Euler method. To get good starting values for the state and control variables we apply an instantaneous control approach. Instantaneous control is closely related to receding horizon control (RHC) or model predictive control (MPC) with finite time horizon \cite{instant1,instant2,instant3,instant4,instant5}, for instance. Here, one has to solve an optimal control problem with purely stationary (i.e. algebraic) constraints in each iteration, where we utilize a trust region method; see, e.g., \cite{conn2000trust}. All other optimal control problems where solved with a version of the limited-memory BFGS algorithm (the so-called L-BFGS-B algorithm) which is able handle simple box constraints; see, e.g., \cite[Chapter 5]{kelley1999iterative}. All codes are implemented in Python. All calculations were performed on a Dell XPS 13 9310 2-in-1 laptop. It contains an Intel i7 11-th generation CPU and 32 GB RAM.

\subsection{Run 1: Small Network}
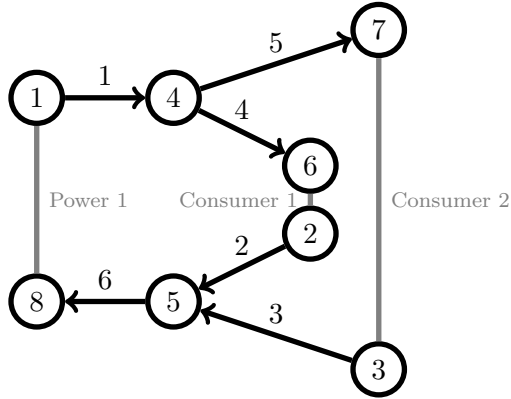
\begin{figure}
\setcounter{figure}{0}
\begin{center}
\begin{tikzpicture}[scale=1.8]
\node[line width=2pt,black] (v1) at (0,1) [circle,draw] {$1$};
\node[line width=2pt,black] (v2) at (2,0) [circle,draw] {$2$};
\node[line width=2pt,black] (v3) at (2.5,-1) [circle,draw] {$3$};
\node[line width=2pt,black] (v4) at (1,1) [circle,draw] {$4$};
\node[line width=2pt,black] (v5) at (1,-0.5) [circle,draw] {$5$};
\node[line width=2pt,black] (v6) at (2,0.5) [circle,draw] {$6$};
\node[line width=2pt,black] (v7) at (2.5,1.5) [circle,draw] {$7$};
\node[line width=2pt,black] (v8) at (0,-0.5) [circle,draw] {$8$};
\draw[->, line width=2pt,] (v1) to node[above] {$1$} (v4);
\draw[->, line width=2pt,] (v2) to node[above] {$2$} (v5);
\draw[->, line width=2pt,] (v3) to node[above] {$3$} (v5);
\draw[->, line width=2pt,] (v4) to node[above] {$4$} (v6);
\draw[->, line width=2pt,] (v4) to node[above] {$5$} (v7);
\draw[->, line width=2pt,] (v5) to node[above] {$6$} (v8);
\draw[-, line width=2pt,gray] (v6) to node[left] {\footnotesize Consumer 1} (v2);
\draw[-, line width=2pt,gray] (v7) to node[right] {\footnotesize Consumer 2} (v3);
\draw[-, line width=2pt,gray] (v1) to node[right] {\footnotesize Power 1} (v8);
\end{tikzpicture}
\end{center}
\caption*{Figure~1 (Run~1): Example of a network with two consumers which are located between the nodes 2 and 6 and 3 and 7. The depot is located between the nodes 1 and 8. The numbers above the arrows correspond to the numbering of the pipes.}
\end{figure}
In this test we consider an example with two consumers, see Figure~1; cf. also \cite{JaekleReichle}. For this network we construct an exact solution with the following parameters
\begin{align*}
k &= -1,    & c_p     &= 2, & d &= 1, & \rho            &= 2, \\
T_\mathsf{ext}&=0, & \lambda &= 2, & L &= 1, & G(\partial_x h) &= 1,
\end{align*}
the consumer demands
\begin{align*}
Q_1(t) = \frac{1}{3}\big(2 \exp(1.5)-1\big) \pi \exp(1+t), \quad Q_2(t) = \frac{1}{6} (2 \exp(3)-1) \pi \exp(1+t)
\end{align*}
and for the velocities, the temperatures and the pressures the functions 
\begin{align*}
&v_1(t)= 3V(t),\, T_1(t,x)= \exp(x-1)T(t),\, 
        p_1(t,0)= 3P(t)+2,\, p_1(t,L) = P(t),  \\
&v_2(t)= 2V(t),\, T_2(t,x) = \frac{1}{2}\exp(\frac{3}{2}x)T(t),\,
        p_2(t,0)= \frac{44}{9}P(t)+4,\, p_2(t,L) = 4P(t)+2,\\
&v_3(t)= V(t),\, T_3(t,x)=\frac{1}{2}\exp(3x)T(t),\, 
        p_3(t,0)=\frac{38}{9}P(t)+4,\,p_3(t,L)=4P(t)+2,\\
&v_4(t)= 2V(t),\, T_4(t,x)= \exp(\frac{3}{2}x)T(t),\, 
        p_4(t,0)= P(t),\, p_4(t,L) =\frac{1}{9}P(t)-2,\\
&v_5(t)= V(t),\, T_5(t,x) = \exp(3x)T(t),\,
        p_5(t,0)= P(t),\, p_5(t,L) = \frac{7}{9}P(t)-2,\\
&v_6(t)= 3V(t), \, T_6(t,x) = \frac{1}{3}(2+\exp(\frac{3}{2}))\exp(\frac{3}{2}+x)T(t)+5, \\
        &p_6(t,0)= 4P(t)+2,\, p_1(t,L) = 2P(t).
\end{align*}
with $V(t) :=1/(6-3t)$, $T(t) :=  \exp(1+t)(2-t)$ and $P(t) = 1/(t-2)^2$.
One can verify that these parameters and equations satisfy the DAE with the underlying differential equations \eqref{eq:endeq}, the coupling conditions \eqref{couplingCond} and the other algebraic equations \eqref{connectNetwork}. With these parameters and functions we can construct control functions which satisfy \eqref{eq:controlFun}. Note that this construction is not unique for $P_w$ and $P_g$, see \eqref{eq:controlFun2} since here the right-hand-side consists of a sum of two controls. Therefore, we will only examine the behavior of $P_p$ here.

For the cost function we set $\omega_i = 0$ for $i=1,2,3$ and for a first test $\alpha$, $\beta$, $\gamma =0.5$ and for a second test $\alpha$, $\beta$, $\gamma = 1$. As desired state $x_d$, the desired algebraic variable $y_d$ and the desired control $u_d$ we use the solution constructed above
\begin{align*}
   x_d(t) &= \begin{pmatrix}
            T_{1,2}(t), \hdots
            ,T_{1,n_1}(t),
            \hdots
            ,T_{6,2}(t),
            \hdots
            ,T_{6,n_6}(t)
         \end{pmatrix}^\top, \\
    y_d(t) &=  \begin{pmatrix}
             v_1(t),
             \hdots
             ,v_6(t),\,
             T_{1,1}(t),
             \hdots
             ,T_{1,6}(t),\,
             p_{1,0}(t),
             \hdots 
             ,p_{6,L_6}(t)
         \end{pmatrix}^\top, \\ 
    u_d(t) &= \begin{pmatrix}
             P_p(t),\,
             P_w(t),\,
             P_g(t)
         \end{pmatrix}^\top
\end{align*}
and set the initial values with $x_0 = x_d(t_0), y_0 = y_d(t_0)$ and $u_0 = u_d(t_0)$. Thus it is clear that the optimal value of the objective function is 0. This value is reached if we use the above constructed controls to calculate the state. 

The time horizon is given by $[t_0, t_f]$ with $t_0 = 0$ and $t_f=1$. Further, we choose the following time grid
\begin{align*}
t_0 < t_1 < \cdots < t_{n_t} = t_f\quad \text{with }n_t\in\mathbb N.
\end{align*}
For the discretization in space we use the method of lines; cf., e.g., \cite{JaekleReichle} for more details. For the discretization in time we use the implicit Euler method with step size $\Delta t = \Delta x \in \{2^{-4}, 2^{-5}, \ldots, 2^{-10} \}$ for both schemes. The behavior of the error in the control for the pressure pumps $P_p$ can be seen in Figure~2.
\begin{figure}
    \setcounter{figure}{1}
	\centering
	\includegraphics[width=65mm,height=50mm]{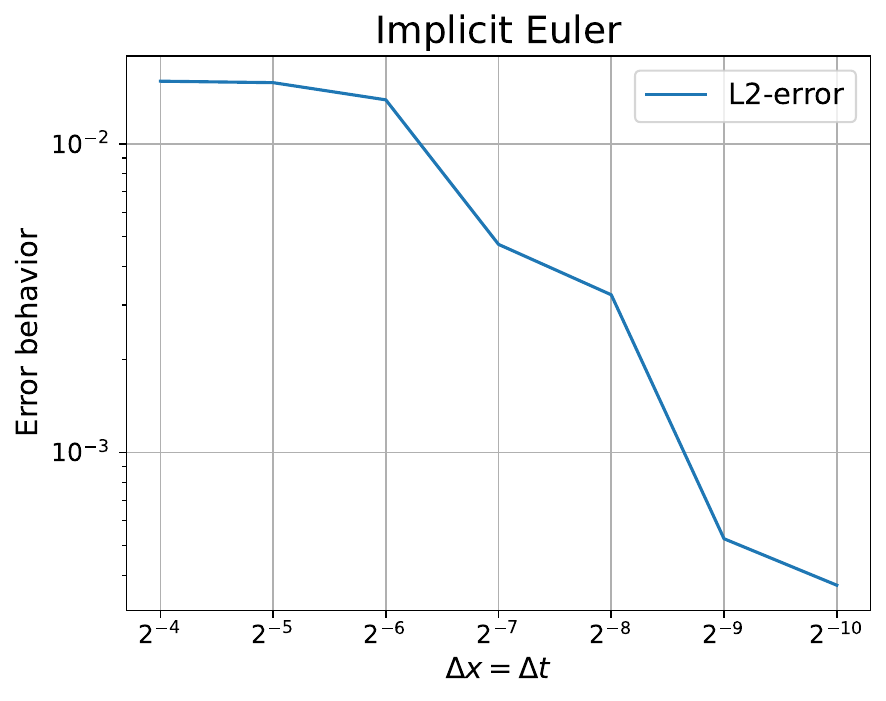}\hspace{5mm}
	\includegraphics[width=65mm,height=50mm]{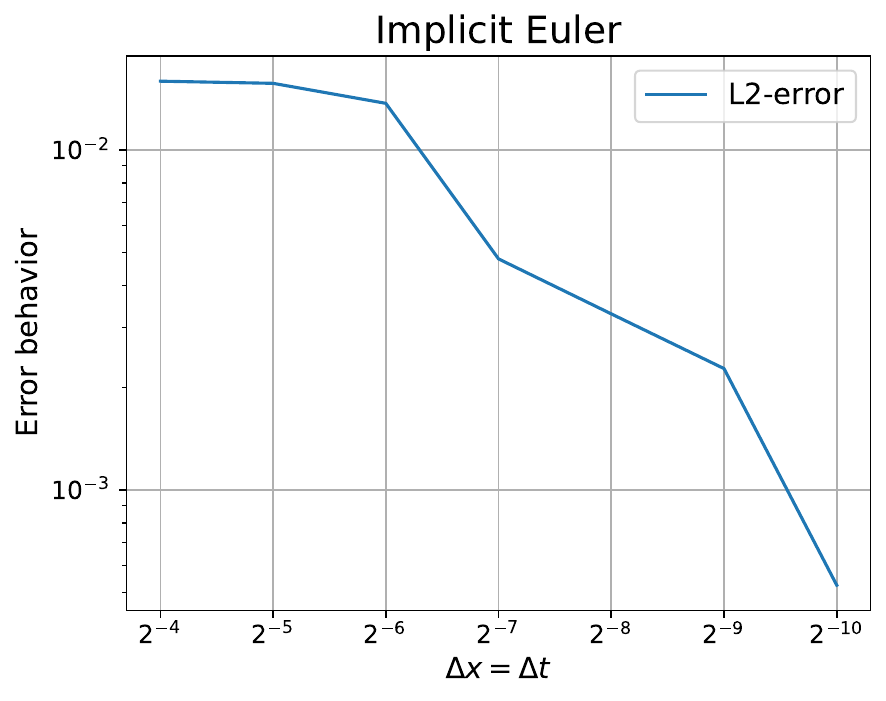}
    \caption*{Figure 2 (Run 1): Discrete $L^2$ error for the control function $P_p$ for $\alpha=\beta=\gamma = 0.5$ (left plot) and $\alpha=\beta=\gamma = 1$ (right plot).}
\end{figure}
As expected (cf., e.g., \cite{Hag00}), we see that the error decreases linearly. Thus, our numerical solution approach performs well.


\subsection{Run 2: Five Consumer}

In this example, we use realistic data to optimize a small network. We want to investigate here how the different parameters in the cost function affect the optimization. The exemplary network is shown in Figure~3.
\begin{figure}
    \setcounter{figure}{2}
	\centering
	\includegraphics[width=75mm,height=50mm]{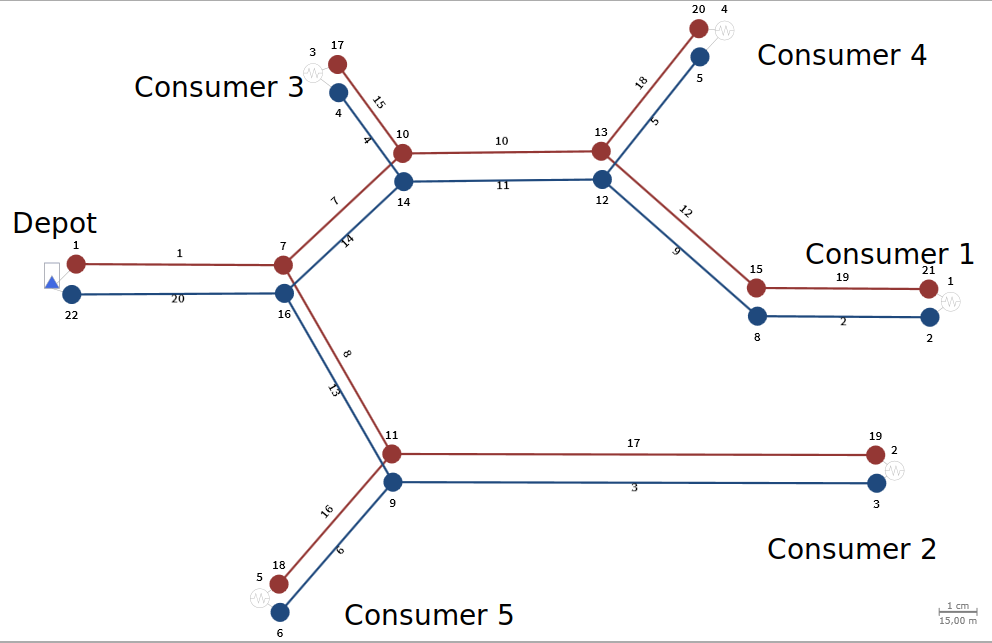}\hspace{5mm}
    \includegraphics[width=55mm,height=50mm]{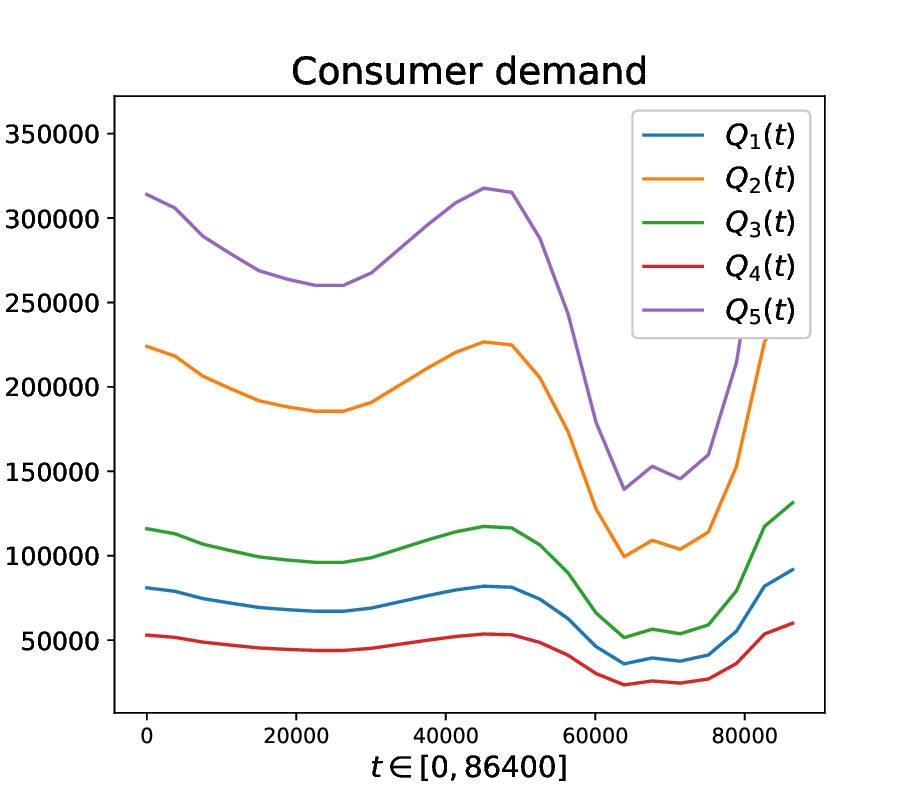}
    \caption*{Figure~3 (Run~2): Network with five consumer, $20$ pipes and $22$ nodes (left plot) and Different consumer demands in the network (right plot)}
\end{figure}
For the instantaneous control approach we consider consumer demand for $24$ hours (i.e., 86400 seconds). However, the optimization is then only for one hour. This can be seen as a first step for a later MPC procedure; cf. \cite{instant2}, for instance. Therefore, the time horizon for the optimization is given by
\begin{align*}
[t_0, t_f] = [0,3600].
\end{align*} 

In the left plot of Figure~3 you can see a numbering of the pipes, nodes and consumers. The red lines and the red circles represent the forward flow in the pipes and nodes, the blue ones the return flow back to the depot. We use for all pipes the parameters
\begin{align*}
    k &= 0.31, & c_p &= 4160, & d &= 0.1071, & \rho &= 960, \hfill\\
    T_\mathsf{ext} &=6, & k_{rough} &= 0.0001, & L &= 100, & \Delta h &= 0.
\end{align*}
The consumption of the five consumers can be seen in the right plot of Figure~3. This reflects the consumption for one day. In the return flow we assume that the water has always $50^\circ\text{C}$. In addition, we set the pressure arriving at the power plant to $2$ bar. Each pipe is discretized into $10$ segments. For the instantaneous control approach we discretize the time interval $[0, 86400]$ into $500$ equidistant pieces and solve therefore $500$ stationary optimization problems. The time interval for the optimization $[0, 3600]$ is discretized into $100$ equidistant pieces.
\begin{table}
    \setcounter{table}{1}
    \centering
        \begin{tabular}{l c c c c} 
            \toprule
            &initial cond&desired variables& optimization&success\\
            \midrule
            $\omega$ equal, reg big & $19$ & $175$ & $121$ & no\\ 
            $\omega$ equal, reg small&$21$&$176$& $220$&yes\\
            $\omega$ different, big equal & $38$ & $170$ &--&no \\
            $\omega$ different, reg small&$41$& $201$ & $172$ &yes\\ [0.1ex] \bottomrule
        \end{tabular}
    \caption*{Table~2 (Run~2): Time in seconds for the different tasks in the optimization procedure. Here `$\omega$ equal' means $\omega_1=\omega_2=\omega_3=10^{-4}$, $\omega$, `$\omega$ different' means  $\omega_1=10^{-4}, \omega_2=10^{-3}, \omega_3=10^{-1}$, `reg small' means $\alpha=\beta=\gamma=10^{-4}$ and `reg big' means $\alpha=10^{-1}, \beta=10^{-1}, \gamma=1$. ' `Success' means here that the termination condition, that the norm of the gradient is small, was reached.}
\end{table}
Table~2 shows comparative values for the CPU time for different optimization situations. We consider on the one hand the simpler case, where the weights of the control (i.e., the $\omega_i$'s) have the same value and on the other hand the more realistic case, where these values are different. In this case, the cost of the pressure pumps is chosen to be the cheapest, while the waste incineration is slightly more expensive and burning gas is the most expensive.

In addition, we distinguish whether we have large or small values for the parameters $\alpha$, $\beta$ and $\gamma$. If these three parameters are chosen to be small, we have less penalization and the solution may be far away from the desired state. Therefore, the solution could be outside the bounds, but the optimization is simplified. If, on the other hand, the parameters $\alpha, \beta$ and $\gamma$ are chosen to be large, the optimization is complicated, which can be noticed in Table~2. In summary, we can state that the most interesting case for the application, different costs for the controls and low values for parameters $\alpha, \beta, \gamma$, gives good results in acceptable time.

In Figure~4 we see the solutions for the instantaneous control approach and the optimal control method.
\begin{figure}
    \setcounter{figure}{4}
    \centering
    \includegraphics[width=65mm,height=50mm]{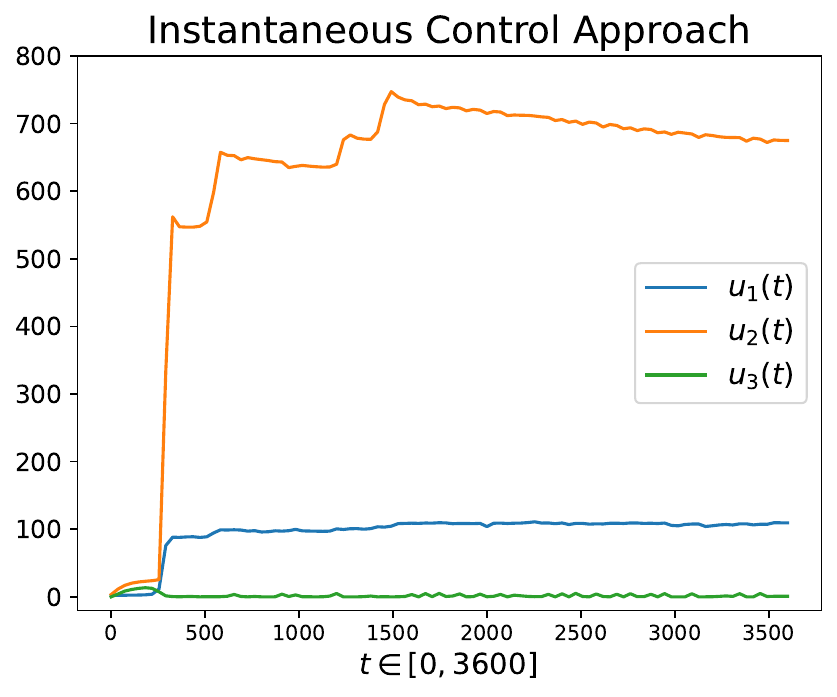}\hspace{5mm}
    \includegraphics[width=65mm,height=50mm]{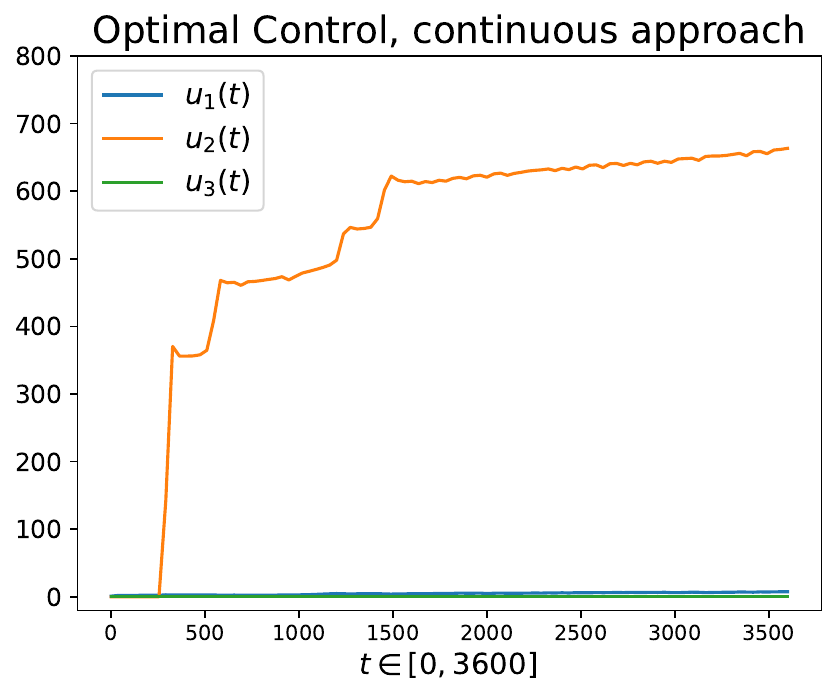}
    \caption*{Figure~4 (Run~2): Control functions for the instantaneous control approach (left plot) and for the optimal control approach (right plot).}
\end{figure}
Here it can be seen very well that the optimal control method finds an even better solution, whereby better means more favorable.

\subsection{Run 3: Part of a real Network}
In the last example, we optimize a subnet of a real network. The data is provided by our project partners (Rechenzentrum f\"ur Versorgungsnetze Wehr GmbH\footnote{See \url{https://www.rzvn.de}}). The layout of the network can be seen in Figure~5.
\begin{figure}
    \setcounter{figure}{5}
    \centering
    \includegraphics[width=120mm,height=70mm]{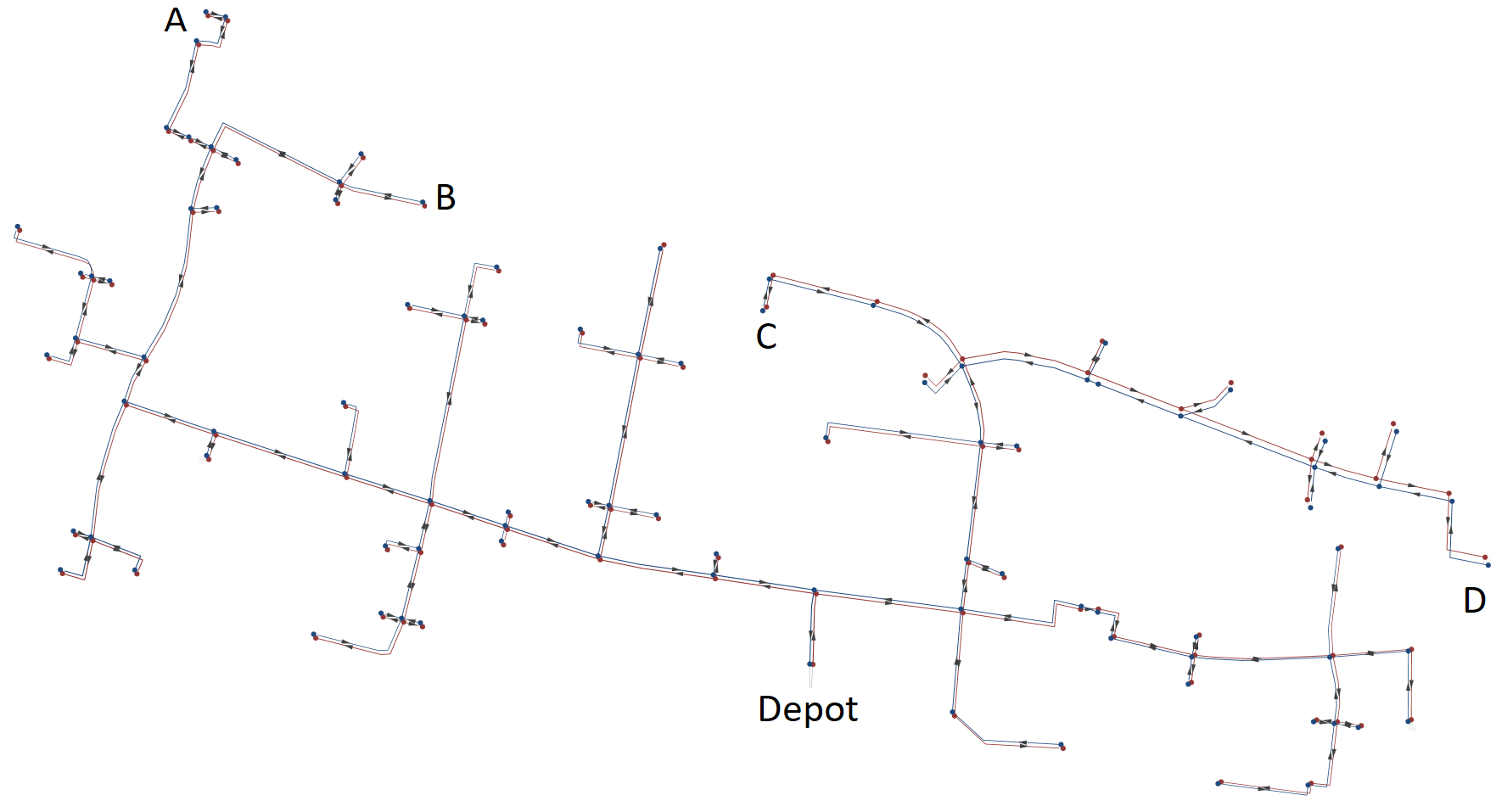}
    \caption*{Figure~5 (Run~3): Part of a real network.}
\end{figure}
The network consists of a total of $193$ nodes and $51$ consumers. Of the $190$ pipes, with a total length of $7988$ meters, $95$ each of the pipes are inflow and return. The total consumption of all consumers for one week can be seen in the left plot of Figure~6. %
\begin{figure}
    \setcounter{figure}{6}
	\centering
	\includegraphics[width=65mm,height=50mm]{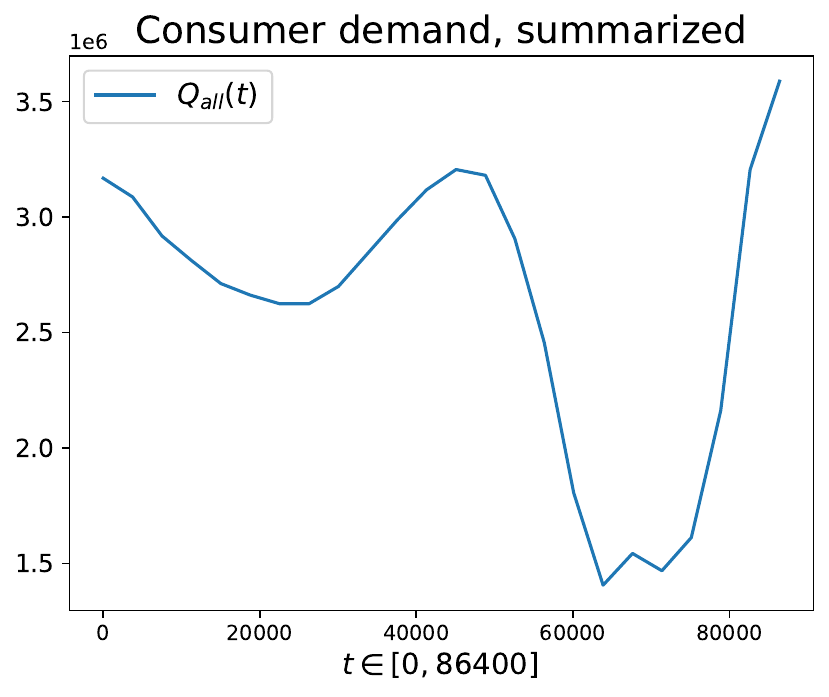}\hspace{5mm}
    \includegraphics[width=65mm,height=50mm]{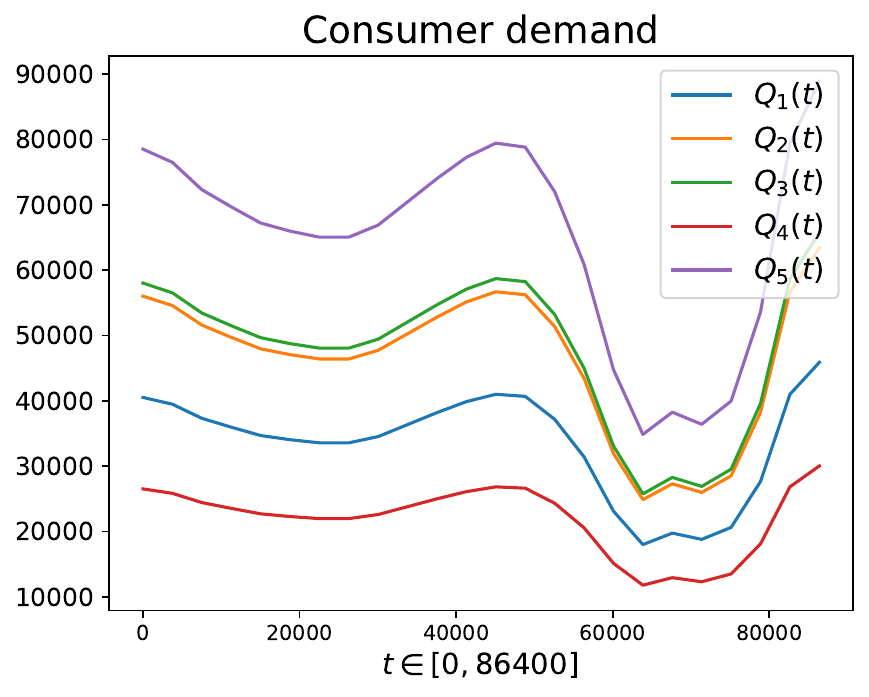}
	\caption*{Figure~6 (Run~3): Total consumption of all consumers for one day (left plot) and consumption for consumers A, B, C and D (right plot).}
	\label{fig:some_consumer}
\end{figure}
The return flow is $50^\circ$C and arrives at the depot with $2$ bar. Since the pipes are buried in the ground, we assume a constant ambient temperature of $6^\circ$C.

We use a coarse space discretization and the pipes are discretized differently. If a pipe is less than $40$ meters long, it is discretized into five sections and all pipes longer than $40$ meters are discretized into twenty sections.

From the results of the previous example, we conclude that we can work with realistic data for the cost function even for a large network. Therefore, we consider here exclusively the case with $\omega_1=10^{-4}$, $\omega_2=10^{-3}$, $\omega_3=10^{-1}$ and $ \alpha=\beta=\gamma=10^{-4}$. As in Run~2 we consider for the instantaneous control approach  a consumer demand for $24$ hours (i.e., 86400 seconds), whereas the optimal control is done only for one hour. For the instantaneous control approach we discretize the time interval $[0, 86400]$ again into $500$ equidistant pieces and for the optimization we discretize the interval $[0, 3600]$ into $100$ equidistant subintervals.

To calculate consistent initial values we need $200$ iterations and $810$ seconds.  The subsequent $500$ stationary optimization problems solved to obtain desired variables take a total of $3250$ seconds. Finally, the L-BFGS-B procedure finds the optimum with a tolerance on the norm of the gradient of $10^{-4}$ in $343$ seconds within $30$ iterations. Note that this is fast enough for a real time optimization in an MPC framework since the optimization for one hour is done in just approximately six minutes. The solution of the optimal control approach is in bounds. 

In Figure~7 we see the solutions for the instantaneous control approach and the optimal control method.
\begin{figure}
    \setcounter{figure}{6}
    \centering
    \includegraphics[width=65mm,height=50mm]{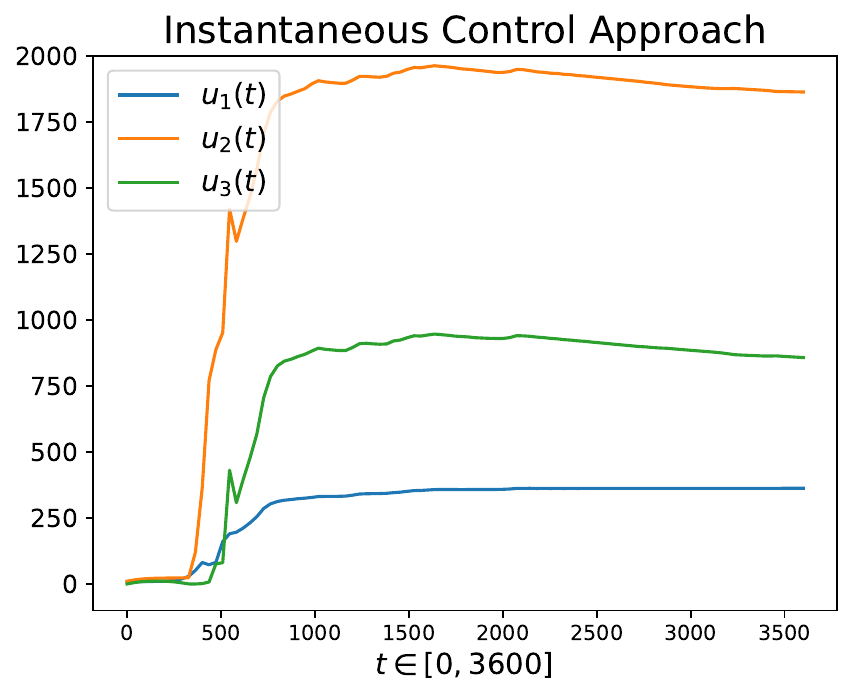}\hspace{5mm}
	\includegraphics[width=65mm,height=50mm]{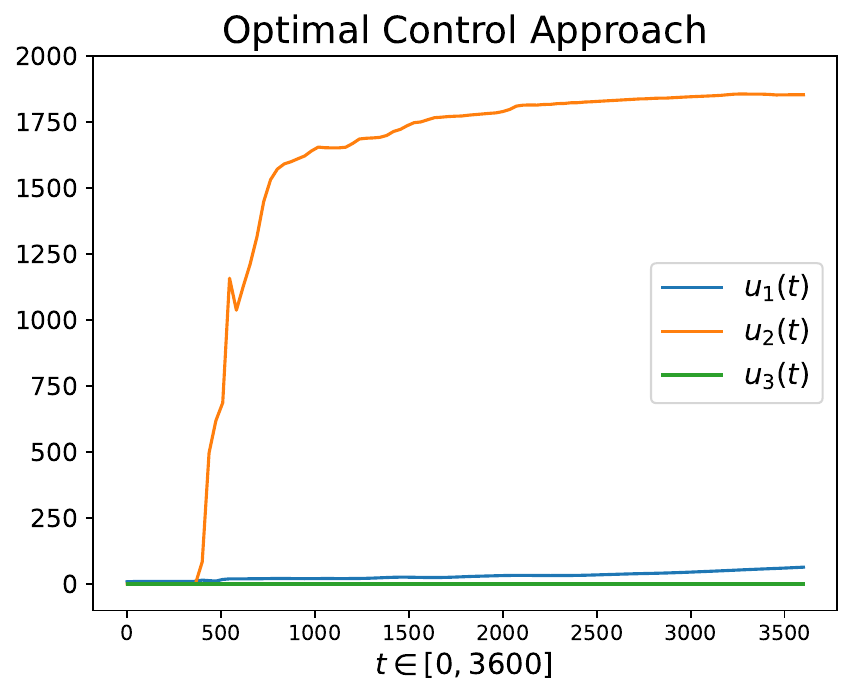}
    \caption*{Figure~7 (Run~3): Control functions for the instantaneous control (left plot) and the optimal control approach (right plot).}
	\label{fig:189Consumer_solution}
\end{figure}
Here it can be seen very well that the optimal control method finds an even better solution, whereby better means more favorable: Comparing the value of the cost without the penalization, e.g.,
\begin{align*}
    \hat J(u):=\int_{t_0}^{t_f} \sum_{i=1}^3 \omega_i(t) u_i(t)\,\mathrm dt
\end{align*}
we observe that the value for the optimal control approach is significantly smaller:
\begin{align*}
    \hat{J}(u_{inst})= 276 093, 07\quad\text{and}\quad\hat{J}(u_{opt}) &= 5 448, 57.
\end{align*}
This again underlines that after finding a suboptimal control with the instantaneous control approach, it is worthwhile to use the optimal control approach. However, the instantaneous control approach is still necessary in our case to get a consistent initial value. Furthermore, we use the solution as desired state and algebraic variable in the cost; compare Remark \ref{rem:desired}.

\section*{Acknowledgements}

We are grateful for Behzad Azmi (Konstanz), Gabriele Ciaramella (Milan), Piet Hensel and Renhang Qin (both from the company RZVN, Düsseldorf; see \texttt{www.rzvn.de}; see \cite{Hen13}) for fruitful discussions and very helpful remarks. Moreover, we thank the company RZVN for providing the data for Run~3.

\section*{Funding}

The authors acknowledge funding by the Federal Ministry of Education and Research (BMBF) for the project \emph{Entwicklung einer Software zur optimalen Auslegung und dynamischen Betriebsweise von bidirektionalen Anergie und W\"armenetzen mit dezentralen Einspeisungen}.

\setcounter{section}{0}
\renewcommand{\theequation}{\thesection.\arabic{equation}}
\renewcommand{\thesection}{\Alph{section}}
    
\section{Proof of Lemma~\ref{lem:bound_existencetime}}
\label{Sec:A}

Due Assumptions~\ref{ass:represent_f}-3) and \ref{ass:representation_yu}-4) there exist $\rho_\circ>0$ and functions $\alpha_1,\alpha_2,\gamma_1,\gamma_2\in L^4(t_0,t_f)$. For $\mathfrak u_\mathsf{ab}=\max\{\|\ua\|_2,\|\ub\|_2\}$ and $u\in\uad$ we have
\begin{align*}
    0&\le\alpha_1\,{\|u(\cdot)\|}_2+\alpha_2\le\alpha_1\mathfrak u_\mathsf{ab}+\alpha_2=:\eta_1\in L^4(t_0,t_f),\\
    0&\le1+\gamma_1\,{\|u(\cdot)\|}_2+\gamma_2\le1+\gamma_1\,\mathfrak u_\mathsf{ab}+\gamma_2=:\eta_2\in L^4(t_0,t_f).
\end{align*}
This implies that $\alpha=\eta_1\eta_2\in L^2(t_0,t_f)$ holds. Hence, we estimate f.a.a. $t\in[t_0,t_f]$ and $\mathsf x,\tilde{\mathsf x}\in B_{\rho_\circ}(x_0)$
\begin{align}
    \label{Estimate-alpha}
    \begin{aligned}
        &{\|f(t,\mathsf x,\mathsf y(t,\mathsf x,u(t)),u(t))-f(t,\tilde{\mathsf x},\mathsf y(t,\tilde{\mathsf x},u(t)),u(t))\|}_2\\
        &\quad\le\big(\alpha_1(t)\,{\|u(t)\|}_2+\alpha_2(t)\big){\|(\mathsf x,\mathsf y(t,\mathsf x,u(t)))-(\tilde{\mathsf x},\mathsf y(t,\tilde{\mathsf x},u(t)))\|}_2\\
        &\quad\le \eta_1(t)\left({\|\mathsf x-\tilde{\mathsf x}\|}_2+{\|\mathsf y(t,\mathsf x,u(t))-\mathsf y(t,\tilde{\mathsf x},u(t))\|}_2\right)\\
        &\quad\le\eta_1(t)\left({\|\mathsf x-\tilde{\mathsf x}\|}_2+{\|\mathsf y_1(t,\mathsf x)-\mathsf y_1(t,\tilde{\mathsf x})\|}_F{\|u(t)\|}_2+{\|\mathsf y_2(t,\mathsf x)-\mathsf y_2(t,\tilde{\mathsf x})\|}_2\right)\\
        &\quad\le\eta_1(t)\left({\|\mathsf x-\tilde{\mathsf x}\|}_2+\gamma_1(t)\,{\|\mathsf x-\tilde{\mathsf x}\|}_2{\|u(t)\|}_2+\gamma_2(t)\,{\|\mathsf x-\tilde{\mathsf x}\|}_2\right)\\
        &\quad\le\eta_1(t)\big( 1+\gamma_1(t)\,{\|u(t)\|}_2+\gamma_2(t)\big){\|\mathsf x-\tilde{\mathsf x}\|}_2\\
        &\quad=\eta_1(t)\eta_2(t){\|\mathsf x-\tilde{\mathsf x}\|}_2=\alpha(t){\|\mathsf x-\tilde{\mathsf x}\|}_2,
    \end{aligned}
\end{align}
i.e., $f$ is locally Lipschitz-continuous on $\dx$ and $\alpha$ does not depend on the chosen $u\in\uad$. Moreover, Assumption~\ref{ass:representation_yu}-5) implies that for our initial value $x_0$ there are $\beta_1,\beta_2\in L^2(t_0,t_f)$ satisfying
\begin{align}
    \label{Estimate-beta}
    \begin{aligned}
        {\|f(t,x_0,\mathsf y(t,x_0,u(t)),u(t))\|}_2&\le\beta_1(t)\,{\|u(t)\|}_2+\beta_2(t)\\
        &\le\beta_1(t)\,\mathfrak u_\mathsf{ab}+\beta_2(t)=:\beta(t),    
    \end{aligned}
\end{align}
where $\beta$ belongs to $L^2(t_0,t_f)$. Now we introduce the non-negative functions 
\begin{align*}
	a(t):= \int_{t_0}^t \alpha(s) \, \mathrm ds\ge0\quad\text{and}\quad
	b(t) := \int_{t_0}^t \beta(s) \,\mathrm ds\ge0
\end{align*}
for $t\in (t_0,t_f]$. To follow the arguments from the proof of Theorem \ref{thm:caratheo_ex} let $t_e\in(t_0,t_f]$ be maximal with
\begin{align*}
	&a(t) \leq a(t_e) < 1 \mbox{ for all } t \in [t_0,t_e],\\
	&a(t)\rho_\circ + b(t) \leq b(t_e)\rho_\circ + b(t_e) \leq \rho_\circ \text{ for all } t\in [t_0,t_e].
\end{align*}
Then with the arguments from the proof of Theorem \ref{thm:caratheo_ex} the claim follows.\hfill\qedsymbol

\section{Proof of Lemma~\ref{lem:controlstate2}}
\label{Sec:B}

Note that the sequence $\{u_k\}_{k\in\mathbb N}$ belongs to $\uad$, which implies $\|u_k\|_\u \leq \mathfrak u_\mathsf{ab}$ with $\mathfrak u_\mathsf{ab}=\max\{\|\ua\|_2,\|\ub\|_2\}$. Therefore, there exists a subsequence -- again denoted by $\{u_k\}_{k\in\mathbb N}$ -- and an element $\tu\in\uad$ with $u_k\rightharpoonup\tu$ for $k\to\infty$. The solution pairs $(x_k,y_k) \in\x\times\y$ are well-defined by Remark~\ref{Rem:controlstate1} for any $k\in\N$. Due to Remark~\ref{rem:a35} and Assumption~\ref{A-3} there is a radius $\rho_\circ>0$ such that $\|x_k(t)-x_0\|_2\le\rho_\circ$ for all $t\in[t_0,t_f]$. From \eqref{Estimate-alpha} and \eqref{Estimate-beta} we infer f.a.a. $t\in[t_0,t_f]$
\begin{align*}
    &{\|f(t,x_k(t),y_k(t),u_k(t))\|}_2^2={\|f(t,x_k(t),\mathsf y(t,x_k(t),u_k(t)),u_k(t))\|}_2^2\\
    &\quad\le 2\,{\|f(t,x_k(t),\mathsf y(t,x_k(t),u_k(t)),u_k(t))-f(t_0,x_0,\mathsf y(t_0,x_k(t_0),u_k(t_0)),u_k(t_0))\|}_2^2\\
    &\qquad+2\,{\|f(t_0,x_0,\mathsf y(t_0,x_k(t_0),u_k(t_0)),u_k(t_0))\|}_2^2\\
    &\quad\leq 2\alpha^2(t)\,{\|x_k(t)-x_0\|}_2^2 + \beta^2(t)\leq 2\alpha^2(t)\rho_\circ^2+2\beta^2(t)=\mu(t).
\end{align*}
Since $\alpha,\beta\in L^2(t_0,t_f)$ holds, we have $\mu\in L^1(t_0,t_f)$. Consequently, we obtain that $f(\cdot,x_k(\cdot),y_k(\cdot),u_k(\cdot)) $ belongs to $L^2(t_0,t_f;\R^\nx)$ for all $k$. Due to $(x_k,y_k)=\mathcal S(u_k)$ it follows that $\dot x_k=f(\cdot,x_k(\cdot),y_k(\cdot),u_k(\cdot))$ a.e. in $(t_0,t_f]$ which implies
\begin{align*}
    {\|\dot x_k\|}_{L^2(t_0,t_f;\R^\nx)}=\bigg(\int_{t_0}^{t_f}{\|f(t,x_k(t),y_k(t),u_k(t))\|}_2^2\,\mathrm dt\bigg)^{1/2}\le{\|\mu\|}_{L^1(t_0,t_f)}^{1/2}<\infty
\end{align*}
for all $k$. Together with \eqref{AprioriEst} we have that $\{x_k\}_{k\in\mathbb N}$ is uniformly bounded in $H^1(t_0,t_f;\R^\nx)$. Thus, there exists a subsequence -- again denoted by $\{x_k\}_{k\in\mathbb N}$ -- and an element $\tx\in H^1(t_0,t_f;\R^\nx)$ such that $x_k\rightharpoonup\tx$ for $k\to\infty$. Moreover, we invoke the compact embedding $H^1(t_0,t_f;\R^\nx) \subset  C([t_0,t_f];\R^\nx)$, i.e., $\{x_k\}_{k\in \N}$ converges (strongly) to $\tx$ in $C([t_0,t_f];\R^\nx)$.\hfill\\
It remains to show, that it holds $\tx = x(\tu)$, which means, that $\tx$ solves \eqref{eq:DAE_Orig} for $u=\tu$. Due to Assumptions~\ref{ass:represent_f}-2) and \ref{ass:representation_yu}-3) functions $f_1(t,\cdot)$, $f_2(t,\cdot)$, $f_3(t,\cdot)$, $\mathsf y_1(t,\cdot)$ and $\mathsf y_2(t,\cdot)$ are continuous on $\dx$ f.a.a. $t\in[t_0,t_f]$. Therefore we have the following pointwise convergences
\begin{align*}
    f_i(t,x_{k_\ell}(t))\to f_i(t,\tx(t))&\quad\text{for }\ell\to\infty\quad\text{f.a.a. }t \in [t_0,t_f],~i=1,2,3,\\
        \mathsf y_j(t,x_{k_\ell}(t)) \to \mathsf y_j(t,\tx(t))&\quad\text{for }\ell\to\infty\quad\text{f.a.a. }t \in [t_0,t_f],~j=1,2.
    \end{align*}
    This implies that f.a.a. $t \in [t_0,t_f]$
    \begin{align}\label{eq:faa_mixed}
        &f_i(t,x_{k_\ell}(t))\mathsf y_j(t,x_{k_\ell}(t)) \to f_i(t,\tx(t))\mathsf y_j(t,\tx(t))&\quad\text{for }\ell\to\infty,
    \end{align}
    for $i=1,2,3$ and $j=1,2$. 
    In the next step, we proof that the convergence \eqref{eq:faa_mixed} is already in $L^2$ by utilizing the dominated convergence theorem \cite[p.~321]{Rud76}. For that purpose we 
    require majorants. Due to Assumption \ref{ass:represent_f}-4) and Assumption \ref{ass:representation_yu}-6) it follows
    \begin{align*}
        {\|f_i(t,\mathsf x)y_j(t,\mathsf x)\|}_2 \leq \iota_i(t)\kappa_j(t)\quad\text{and}\quad {\|f_i(t,\mathsf x)\|}_2 \leq \iota_i(t)
    \end{align*}
    f.a.a. $t \in [t_0,t_f]$, for all $\mathsf x \in B_{\rho_\circ}(x_0)$ and for $i = 1,2,3$, $j = 1,2$. Hölder's inequality implies that $\iota_i\kappa_j \in L^2(t_0,t_f)$ for $i = 1,2,3$ and $j = 1,2$. 
    Due to the dominated convergence theorem we have
    \begin{align}
        \label{Theorem5.12}
        \begin{aligned}
            \int_{t_0}^{t_f}{\|f_i(t,x_{k_\ell}(t))\mathsf y_j(t,x_{k_\ell}(t))-f_i(t,\tx(t))\mathsf y_j(t,\tx(t))\|}_2^2\,\mathrm dt&\stackrel{\ell\to\infty}{\longrightarrow}0,\\
            \int_{t_0}^{t_f}{\|f_i(t,x_{k_\ell}(t))-f_i(t,\tx(t))\|}_2^2\,\mathrm dt&\stackrel{\ell\to\infty}{\longrightarrow}0,
        \end{aligned}
    \end{align}
    for $i=1,2,3$ and $j=1,2$.
    For all $t \in [t_0,t_f]$ we also infer
    \begin{align}
        \label{esti-1}
        \begin{aligned}
            \int_{t_0}^t f_1(s,x_{k_\ell}(s))\mathsf y_2(s,x_{k_\ell}(s)) \, \mathrm d s&\stackrel{\ell\to\infty}{\longrightarrow} \int_{t_0}^t f_1(s,\tx(s))\mathsf y_2(s,\tx(s))\,\mathsf d s,\\
		    \int_{t_0}^t f_3(s,x_{k_\ell}(s))\,\mathrm d s&\stackrel{\ell\to\infty}{\longrightarrow}\int_{t_0}^t f_3(s,\tx(s))\,\mathrm ds.
	   \end{aligned}
    \end{align}
    Due to \eqref{Theorem5.12} and $u_{k_\ell}\rightharpoonup\tu$ for $\ell\to\infty$ in $\u$ we conclude from \cite[Theorem 5.12-4 (c)]{Ciarlet13} that for all $t \in [t_0,t_f]$
    \begin{align}
        \label{esti-2}
        \begin{aligned}
            \int_{t_0}^t f_1(s,x_{k_\ell}(s))\mathsf y_1(s,x_{k_\ell}(s))u_{k_\ell}(s) \,\mathsf d s 
            &\stackrel{\ell\to\infty}{\longrightarrow}\int_{t_0}^t f_1(s,\tx(s))\mathsf y_1(s,\tx(s))\tu(s) \,\mathrm d s,\\
		  \int_{t_0}^t f_2(s,x_{k_\ell}(s))u_{k_\ell}(s) \,\mathrm d s&\stackrel{\ell\to\infty}{\longrightarrow}\int_{t_0}^t f_2(s,\tx(s))\tu(s) \, \mathrm d s.
	   \end{aligned}
    \end{align}
    Now we get that $\tx$ satisfies the DAE system. In fact, using Lemma~\ref{lem:ac_derivative}, $(x_{k_\ell},y_{k_\ell})=\mathcal S(u_{k_\ell})$, \eqref{esti-1} and \eqref{esti-2} we obtain for all $t \in [t_0,t_f]$ 
    \begin{align*}
		\tx(t)&=\lim_{\ell\to\infty}x_{k_\ell}(t)=\lim_{k\to\infty}\bigg(x_{k_\ell}(t_0)+\int_{t_0}^t f(s,x_{k_\ell}(s),y_{k_\ell}(s),u_{k_\ell}(s))\,\mathrm d s\bigg)\\ 
        &=\lim_{k\to\infty}\bigg(x_{k_\ell}(t_0)+\int_{t_0}^t f(s,x_{k_\ell}(s),\mathsf y(s,x_{k_\ell}(s),u_{k_\ell}(s)),u_{k_\ell}(s))\,\mathrm d s\bigg)\\
		&=\lim_{k\to\infty}\bigg(x_{k_\ell}(t_0)+\int_{t_0}^t f_1(s,x_{k_\ell}(s))\mathsf y(s,x_{k_\ell}(s),u_{k_\ell}(s))+f_2(s,x_{k_\ell}(s))u_{k_\ell}(s)\,\mathrm d s\\ 
        &\hspace{15mm}+\int_{t_0}^tf_3(s,x_{k_\ell}(s))\,\mathrm d s\bigg)\\
		&=\lim_{k\to\infty}\bigg(x_{k_\ell}(t_0)+\int_{t_0}^t f_1(s,x_{k_\ell}(s))\mathsf y_1(s,x_{k_\ell}(s))u_{k_\ell}(s)\,\mathrm d s\\ 
        &\hspace{15mm}+\int_{t_0}^tf_1(s,x_{k_\ell}(s))\mathsf y_2(s,x_{k_\ell}(s))+f_2(s,x_{k_\ell}(s))u_{k_\ell}(s)+ f_3(s,x_{k_\ell}(s))\,\mathrm d s\bigg)\\
        &=\tx(t_0)+\int_{t_0}^t f(s,\tx(s),\ty(s),\tu(s))\,\mathrm ds.
	\end{align*} 
	Consequently, $\tx$ is absolutely continuous, has derivative $\dot{\tx} = f(\cdot,\tx(\cdot),\ty(\cdot),\tu(\cdot))$ in $(t_0,t_f]$ and satisfies $\tx(t_0) = x_0$. 
	Therefore we have $\tx = x(\tu)$, which completes our proof.\hfill\qedsymbol

\nocite{*}
\bibliographystyle{plain}
\bibliography{references.bib}

\begin{thebibliography}{10}

\bibitem{Banda06}
Mapundi Banda, Michael Herty, and Axel Klar.
\newblock Coupling conditions for gas networks governed by the isothermal euler
  equations.
\newblock {\em NHM}, 1:295--314, 06 2006.

\bibitem{betts2001practical}
John~T. Betts.
\newblock {\em Practical methods for optimal control using nonlinear
  programming, ser}, volume~36 of {\em Advances in Design and Control}.
\newblock Society for Industrial and Applied Mathematics (SIAM), Philadelphia,
  PA, thrid edition, 2020.

\bibitem{instant5}
Thomas~R. Bewley, Parviz Moin, and Roger Temam.
\newblock {DNS}-based predictive control of turbulence: an optimal benchmark
  for feedback algorithms.
\newblock {\em Journal of Fluid Mechanics}, 447:179--225, 2001.

\bibitem{BGV16}
Chiara Bordin, Angelo Giordini, and Daniele Vigo.
\newblock An optimization approach for district heating strategic network
  design.
\newblock {\em European Journal of Operational Research}, 252:296--307, 2016.

\bibitem{Borsche19}
Raul Borsche, Matthias Eimer, and Norbert Siedow.
\newblock A local time stepping method for thermal energy transport in district
  heating networks.
\newblock {\em Applied Mathematics and Computation}, 353:215--229, 07 2019.

\bibitem{bressan2007introduction}
Alberto Bressan and Benedetto Piccoli.
\newblock {\em Introduction to the mathematical theory of control}, volume~1.
\newblock American institute of mathematical sciences Springfield, 2007.

\bibitem{gabriele}
Gabriele Ciaramella.
\newblock Optimal control of ordinary differential equations.
\newblock {\em Lecture Notes}, 2018.

\bibitem{Ciarlet13}
Philippe~G. Ciarlet.
\newblock {\em Linear and Nonlinear Functional Analysis with Applications}.
\newblock Society for Industrial and Applied Mathematics, Philadelphia, USA,
  2013.

\bibitem{colella2012numerical}
Francesco Colella, Adriano Sciacovelli, and Vittorio Verda.
\newblock Numerical analysis of a medium scale latent energy storage unit for
  district heating systems.
\newblock {\em Energy}, 45(1):397--406, 2012.

\bibitem{Colombo08}
Rinaldo~M. Colombo and Mauro Garavello.
\newblock On the cauchy problem for the p-system at a junction.
\newblock {\em SIAM Journal on Mathematical Analysis}, 39(5):1456--1471, 2008.

\bibitem{conn2000trust}
Andrew~R. Conn, Nicholas~I.M. Gould, and Philippe~L. Toint.
\newblock {\em Trust region methods}.
\newblock SIAM, 2000.

\bibitem{DMRS23}
Hannes D\"anschel, Volker Mehrmann, Marius Roland, and Martin Schmidt.
\newblock Adaptive nonlinear optimization of district heating networks based on
  model and discretization catalogs.
\newblock {\em SeMA}, 2023.

\bibitem{Domschke15}
Pia Domschke, Oliver Kolb, and Jens Lang.
\newblock Adjoint-based error control for the simulation and optimization of
  gas and water supply networks.
\newblock {\em Applied Mathematics and Computation}, 259:1003--1018, 05 2015.

\bibitem{dorfner2014large}
Johannes Dorfner and Thomas Hamacher.
\newblock Large-scale district heating network optimization.
\newblock {\em IEEE Transactions on Smart Grid}, 5(4):1884--1891, 2014.

\bibitem{Eva10}
Lawrence~C. Evans.
\newblock {\em Partial Differential Equations}.
\newblock Graduate Studies in Mathematics. American Mathematical Society,
  Providence, Rhode Island, 2010.

\bibitem{filippov1962certain}
Aleksei~F. Filippov.
\newblock On certain questions in the theory of optimal control.
\newblock {\em Journal of the Society for Industrial and Applied Mathematics,
  Series A: Control}, 1(1):76--84, 1962.

\bibitem{gerdts2006local}
Matthias Gerdts.
\newblock Local minimum principle for optimal control problems subject to
  differential-algebraic equations of index two.
\newblock {\em Journal of Optimization Theory and Applications}, 130:443--462,
  2006.

\bibitem{gerdts2011optimal}
Matthias Gerdts.
\newblock {\em Optimal control of ODEs and DAEs}.
\newblock Walter de Gruyter, 2011.

\bibitem{gerhardt06}
Claus Gerhardt.
\newblock {\em Analysis II}, volume~1.
\newblock International Press, 2006.

\bibitem{instant2}
Lars Gr{\"u}ne and J{\"u}rgen Pannek.
\newblock {\em Nonlinear Model Predictive Control}.
\newblock Communications and Control Engineering (CSE). Springer, 2017.

\bibitem{Hag00}
William~W. Hager.
\newblock Runge-{K}utta methods in optimal control and the transformed adjoint
  system.
\newblock {\em Numerische Mathematik}, 87:247–282, 2000.

\bibitem{hart2011pyomo}
William~E Hart, Jean-Paul Watson, and David~L Woodruff.
\newblock Pyomo: modeling and solving mathematical programs in python.
\newblock {\em Mathematical Programming Computation}, 3:219--260, 2011.

\bibitem{hauschild2020port}
Sarah-Alexa Hauschild, Nicole Marheineke, Volker Mehrmann, Jan Mohring,
  Arbi~Moses Badlyan, Markus Rein, and Martin Schmidt.
\newblock Port-hamiltonian modeling of district heating networks.
\newblock In {\em Progress in Differential-Algebraic Equations II}, pages
  333--355, Cham, 2020. Springer International Publishing.

\bibitem{Hen13}
Piet Hensel.
\newblock {\em Optimierung des Ausbaus von Nah- und Fernw{\"a}rmenetzen --
  unter Ber{\"u}ck\-sichtigung eines bestehenden Gasnetzes}.
\newblock PhD thesis, Universit{\"a}t Paderborn, Germany, 2013.

\bibitem{instant4}
Michael Hinze.
\newblock Instantaneous closed loop control of the {N}avier-{S}tokes system.
\newblock {\em SIAM journal on control and optimization}, 44(2):564--583, 2005.

\bibitem{instant1}
Michael Hinze and Stefan Volkwein.
\newblock Analysis of instantaneous control for the {B}urgers equation.
\newblock {\em Nonlinear Analysis: Theory, Methods \& Applications},
  50(1):1--26, 2002.

\bibitem{instant3}
Kazufumi Ito and Karl Kunisch.
\newblock Receding horizon optimal control for infinite dimensional systems.
\newblock {\em ESAIM: control, optimisation and calculus of variations},
  8:741--760, 2002.

\bibitem{JaekleReichle}
Christian Jäkle, Lena Reichle, and Stefan Volkwein.
\newblock Modelling and simulation of district heating networks.
\newblock 2023.

\bibitem{kelley1999iterative}
Carl~T. Kelley.
\newblock {\em Iterative methods for optimization}.
\newblock SIAM, 1999.

\bibitem{Koecher00}
Ralf K\"ocher.
\newblock {\em Beitrag zur Berechnung und Auslegung von {F}ernw\"armenetzen}.
\newblock PhD thesis, Technische Universit{\"a}t Berlin,
  Universit{\"a}tsbibliothek (Diss.-Stelle), 01 2000.

\bibitem{Krug19}
Richard Krug, Volker Mehrmann, and Martin Schmidt.
\newblock Nonlinear optimization of district heating networks.
\newblock {\em Optimization and Engineering}, 22:783–819, 2021.

\bibitem{KM06}
Peter Kunkel and Volker Mehrmann.
\newblock {\em Differential-Algebraic Equations: Analysis and Numerical
  Solution}.
\newblock European Mathematical Society Press, Berlin, 2016.

\bibitem{Leo17}
Giovanni Leoni.
\newblock {\em A First Course in Sobolev Spaces}, volume 181 of {\em Graduate
  Studies in Mathematics}.
\newblock American Mathematical Society, 2 edition, 2017.

\bibitem{martens2019convergence}
Bj{\"o}rn Martens and Matthias Gerdts.
\newblock Convergence analysis of the implicit euler-discretization and
  sufficient conditions for optimal control problems subject to index-one
  differential-algebraic equations.
\newblock {\em Set-Valued and Variational Analysis}, 27:405--431, 2019.

\bibitem{martens2021convergence}
Bj{\"o}rn Martens and Matthias Gerdts.
\newblock Convergence analysis for approximations of optimal control problems
  subject to higher index differential-algebraic equations and pure state
  constraints.
\newblock {\em SIAM Journal on Control and Optimization}, 59(3):1903--1926,
  2021.

\bibitem{OR00}
James~M. Ortega and Werner~C. Rheinboldt.
\newblock {\em Iterative Solution of Nonlinear Equations in Several Variables}.
\newblock Classics in Applied Mathematics. Society for Industrial and Applied
  Mathematics, 2000.

\bibitem{pirouti2013energy}
Marouf Pirouti, Audrius Bagdanavicius, Janaka Ekanayake, Jianzhong Wu, and Nick
  Jenkins.
\newblock Energy consumption and economic analyses of a district heating
  network.
\newblock {\em Energy}, 57:149--159, 2013.

\bibitem{rein2018parametric}
Markus Rein, Jan Mohring, Tobias Damm, and Axel Klar.
\newblock Parametric model order reduction for district heating networks.
\newblock {\em PAMM}, 18(1):e201800192, 2018.

\bibitem{rein2020optimal}
Markus Rein, Jan Mohring, Tobias Damm, and Axel Klar.
\newblock Optimal control of district heating networks using a reduced order
  model.
\newblock {\em Optimal Control Applications and Methods}, 41(4):1352--1370,
  2020.

\bibitem{RMDK21}
Markus Rein, Jan Mohring, Tobias Damm, and Axel Klar.
\newblock Model order reduction of hyperbolic systems focusing on district
  heating networks.
\newblock {\em Journal of the Franklin Institute}, 358:7674--7697, 2021.

\bibitem{rezaie2012district}
Behnaz Rezaie and Marc~A Rosen.
\newblock District heating and cooling: Review of technology and potential
  enhancements.
\newblock {\em Applied energy}, 93:2--10, 2012.

\bibitem{RS20}
Marius Roland and Martin Schmidt.
\newblock Mixed-integer nonlinear optimization for district heating network
  expansion.
\newblock {\em at -- Automatisierungstechnik}, 68:985--1000, 2020.

\bibitem{roubivcek2002optimal}
Tom{\'{a}}{\v{s}} Roubi{\v{c}}ek and Michael Val{\'{a}}{\v{s}}ek.
\newblock Optimal control of causal differential--algebraic systems.
\newblock {\em Journal of Mathematical Analysis and Applications},
  269(2):616--641, 2002.

\bibitem{roxin1962existence}
Emilio Roxin.
\newblock The existence of optimal controls.
\newblock {\em Michigan Mathematical Journal}, 9(2):109--119, 1962.

\bibitem{Rud76}
Walter Rudin.
\newblock {\em Principles of Mathematical Analysis}.
\newblock International series in pure and applied mathematics. McGraw-Hill,
  1976.

\bibitem{schweiger2017district}
Gerald Schweiger, Per-Ola Larsson, Fredrik Magnusson, Patrick Lauenburg, and
  St{\'e}phane Velut.
\newblock District heating and cooling systems--framework for modelica-based
  simulation and dynamic optimization.
\newblock {\em Energy}, 137:566--578, 2017.

\bibitem{sontag1998}
Eduardo~D. Sontag.
\newblock {\em Mathematical Control Theory -- Deterministic Finite Dimensional
  Systems}, volume~2.
\newblock Springer New York, NY, 1998.

\bibitem{verda2011primary}
Vittorio Verda and Francesco Colella.
\newblock Primary energy savings through thermal storage in district heating
  networks.
\newblock {\em Energy}, 36(7):4278--4286, 2011.

\end{thebibliography}

\end{document}